\documentclass[11pt]{article}

\usepackage{hyperref}
\usepackage{array} 
\usepackage{amsfonts}
\usepackage{amsmath}
\usepackage{amssymb}
\usepackage{graphicx}
\usepackage{color}

\newtheorem{thm}{Theorem}[section]
\newtheorem{cor}[thm]{Corollary}
\newtheorem{lem}[thm]{Lemma}
\newtheorem{prop}[thm]{Proposition}
\newtheorem{defn}[thm]{Definition}

\numberwithin{equation}{section}

\newcommand{\dx}{\,{\rm d}x}
\newcommand{\dy}{\,{\rm d}y}

\newcommand{\dt}{\,{\rm d}t}
\newcommand{\rd}{{\rm d}}

\def\LL{\mathrm{L}} 
\def\supp{\mathrm{supp}} 

\newcommand{\A}{\mathcal{L}}
\newcommand{\AM}{\mathcal{L}^{\frac{1}{2}}}
\newcommand{\AI}{\mathcal{L}^{-1}}
\newcommand{\AIM}{\mathcal{L}^{-\frac{1}{2}}}
\newcommand{\RR}{\mathbb{R}}

\newcommand{\B}{\mathcal{B}}

\def\ee{\mathrm{e}} 
\def\dist{\mathrm{dist}} 
\def\diam{\mathrm{diam}} 

\def\qed{\,\unskip\kern 6pt \penalty 500
\raise -2pt\hbox{\vrule \vbox to8pt{\hrule width 6pt
\vfill\hrule}\vrule}\par}
\definecolor{darkblue}{rgb}{0.05, .05, .65}
\definecolor{darkgreen}{rgb}{0.1, .65, .1}
\definecolor{darkred}{rgb}{0.8,0,0}

\begin{document}
\title{\textbf{A Priori Estimates  for Fractional Nonlinear \\ Degenerate Diffusion Equations \\ on bounded domains}\\[7mm]}

\author{\Large Matteo Bonforte$^{\,a,\,b}$ 
~and~ Juan Luis V\'azquez$^{\,a,\,c}$\\} 
\date{} 

\maketitle

\

\begin{abstract}
We investigate  quantitative properties of the nonnegative solutions
$u(t,x)\ge 0$ to the nonlinear fractional diffusion equation, $\partial_t u
+ \A (u^m)=0$, posed in a bounded domain, $x\in\Omega\subset \RR^N$ with $m>1$ for $t>0$.   As $\A$ we use one of the most common definitions of the fractional Laplacian $(-\Delta)^s$, $0<s<1$,  in a bounded domain with zero Dirichlet boundary conditions. We consider a general class of very weak solutions of the equation, and obtain a priori estimates in the form of smoothing effects, absolute upper bounds, lower bounds, and Harnack inequalities. We also investigate the boundary behaviour and we obtain sharp estimates from above and below.  The standard Laplacian case $s=1$ or the linear case $m=1$ are recovered as limits. The method is quite general, suitable to be applied to a number of similar problems.

\end{abstract}
\vspace{2cm}

\noindent {\rule{\textwidth}{0.5pt}}\\[.2cm]

\noindent {\rule{\textwidth}{0.5pt}}
\begin{itemize}
\item[(a)] Departamento de Matem\'{a}ticas, Universidad
Aut\'{o}noma de Madrid,\\ Campus de Cantoblanco, 28049 Madrid, Spain
\item[(b)] e-mail address:~\texttt{matteo.bonforte@uam.es }\\ web-page:~\texttt{http://www.uam.es/matteo.bonforte}
\item[(c)] e-mail address:~\texttt{juanluis.vazquez@uam.es }\\ web-page:~\texttt{http://www.uam.es/juanluis.vazquez}
\end{itemize}


\vskip .3cm


\newpage
\tableofcontents
\normalsize
\newpage

\section{Introduction}

A cornerstone in the modern theory of nonlinear partial differential equations is the derivation of suitable a priori estimates. Such estimates are used to ensure the boundedness or compactness of the approximate families on which the existence theory is based, and also in deriving the regularity estimates that ensure that generalized solutions are actually smooth in some sense. In this paper we address the question of  obtaining a priori estimates for the  solutions of the Nonlinear  Fractional Diffusion Equation (NFDE):
\begin{equation}\label{FPME.equation}
\partial_t u +\A(u^m)=0  \qquad\mbox{posed in }Q=(0,\infty)\times \Omega\,,
\end{equation}
where $\Omega\subset \RR^N$ is a bounded domain with smooth boundary. We take $m>1$,  $0<s\le 1$, $N\ge 1$ and $u\ge 0$. The linear operator $\A$ is a fractional power of the Laplacian subject to suitable Dirichlet boundary conditions. It is known that different versions are acceptable to materialize such concept of fractional operator. Their definition  and main properties are recalled for the reader's convenience in Section \ref{ssect.Def.Fract.Lapl}. For definiteness, in this paper we will consider one of such definitions, the so-called spectral fractional Laplacian, and derive a whole set of estimates for the nonnegative solutions of the problem, but the arguments are of a general nature and are suitable to be applied to other choices of diffusion operator.

\medskip

\noindent {\sc The main results.}  For a quite general class of nonnegative weak solutions to the above problem, we will derive in Theorems \ref{thm.Upper.PME} and \ref{thm.Upper.PME.II}  absolute upper estimates, valid up to the boundary, of the form
\begin{equation}\label{Intro.Abs.Bdds}
u(t,x) \le K_2\, \dfrac{\Phi_1(x)^{\frac{1}{m}}}{t^{\frac{1}{m-1}}}\,,\qquad\forall t>0\,,\;\forall x\in \Omega\,.
\end{equation}
In particular, we observe that the boundary behaviour is dictated by $\Phi_1$\,, the first positive eigenfunction of $\A$\,, which behaves like (a power of) the distance to the boundary at least when the domain is smooth enough. In Theorem \ref{thm.Upper.1.PME} we will also prove standard and weighted instantaneous smoothing effects of the type
\begin{equation}\label{Intro.smoothing}
\sup_{x\in \Omega}u(t,x)\le \dfrac{K_4}{t^{N\vartheta_{1,1}}}\left(\int_\Omega u(t,x)\Phi_1(x)\dx\right)^{2s\vartheta_{1,1}}
\le\dfrac{K_4}{t^{N\vartheta_{1,1}}}\left(\int_\Omega u_0\Phi_1(x)\dx\right)^{2s\vartheta_{1,1}}\,,\qquad\mbox{$\forall t>0$}\,,
\end{equation}
where $\vartheta_{1,1}=1/(2s+(N+1)(m-1))$\,. This is sharper than \eqref{Intro.Abs.Bdds} only for small times.  As a consequence of the above upper estimates, we derive a number of useful weighted estimates in Section \ref{Sect.Backward}\,, and we also obtain backward in time smoothing effect of the form
\begin{equation}\label{Intro.smoothing.back}
\| u(t)\|_{\LL^\infty(\Omega)} \le \frac{K_4}{t^{(d+1)\vartheta_{1,1}}}\left(1\vee \frac{h}{t}\right)^{\frac{2s\vartheta_{1,1}}{m-1}}\|u(t+h)\|_{\LL^1_{\Phi_1}(\Omega)}^{2s\vartheta_{1,1}}\,,\qquad\mbox{$\forall t,h>0$}\,,
\end{equation}
which is a bit surprising and has not been observed before to our knowledge.
Section \ref{sect.positivity} is devoted to the lower estimates, the main result being Theorem \ref{thm.Lower.PME}:
\begin{equation}\label{Intro.lower}
u(t_0,x_0) \ge L_1\, \dfrac{\Phi_1(x_0)^{\frac{1}{m}}}{t^{\frac{1}{m-1}}}\qquad\mbox{$\forall t\ge t_* > 0\,,\;\forall x_0\in \Omega$}\,,
\end{equation}
where  the waiting time $t_*$  has the explicit form
\[
t_*= L_0\left(\int_{\Omega}u_0\Phi_1\dx\right)^{-(m-1)}\,.
\]
In Section \ref{sect.Haranck} we observe that the above estimates combine into global Harnack inequalities in the form given in Theorem \ref{thm.GHP.PME}\,, namely for all $t\ge t_*$
\begin{equation}\label{Intro.GHP}
H_0\,\frac{\Phi_1(x_0)^{\frac{1}{m}}}{t^{\frac{1}{m-1}}} \le \,u(t,x_0)\le H_1\, \frac{\Phi_1(x_0)^{\frac{1}{m}}}{t^{\frac{1}{m-1}}}\,.
\end{equation}
We provide as a corollary the local Harnack inequalities of elliptic type of Theorem \ref{thm.Harnack.Local}: for all $t\ge t_*$
\begin{equation}\label{Intro.Harnack}
\sup_{x\in B_R(x_0)}u(t,x)\le H_2\,\inf_{x\in B_R(x_0)}u(t,x).
\end{equation}
The latter inequalities imply the more standard forward Harnack inequalities of \cite{D88,DGVbook}, and to our knowledge have not been observed before, even in the case $s=1$\,.

All the above constants $K_i,L_i,H_i>0$ are universal, in the sense that may depend only on $N, m, s$ and $\Omega$\,, but not on $u$. They also have an almost explicit expression, usually given in the proof. Actually, we have tried to obtain quantitative versions of the estimates with indication of the dependence of the relevant constants. In some cases the estimates are absolute, in the sense that they are valid independently of the (norm of the) initial data.

In Section \ref{sect.Elliptic.Prob} we show how our method applies also to weak solutions of elliptic equations of the form
\[
\A W=\lambda W^p\qquad\qquad\mbox{with $p=\frac{1}{m}<1$\,,}
\]
and with homogeneous Dirichlet boundary conditions. In Theorem \ref{Thm.Elliptic.Harnack} we obtain sharp upper bounds and lower bounds for $V=W^p$\,, namely
\[
h_0\|W\|^p_{\LL^p_{\Phi_1}(\Omega)}\Phi_1(x_0)\le W(x_0)\le h_1\Phi_1(x_0)\,.
\]
Note that the upper bound in the formula does not hold for $p=1$.

\medskip

\noindent {\sc New method.} It is worth making a comment about the novel tools used in the proofs. Departing from previous works, we exploit the functional properties of the linear operator as much as possible. More precisely, we use the Green function of the fractional operator  even in the definition of solution, and we make essential use of estimates on its behaviour in the proofs. A key ingredient is thus the knowledge of good estimates for the Green function, that we state in Section \ref{ssect.Funct.Setup}; the new ones are proved in Appendix~\ref{SSection.Green}. Such a functional approach allows us to avoid more standard methods, such as Moser's iteration, the De Giorgi method, or sophisticated Aleksandrov's moving plane methods used for instance in \cite{BV, BV-ADV, BV2012}. A careful inspection of the proofs shows that the present method would allow to treat a quite wide class of linear operators, an issue that we shall discuss in a forthcoming paper \cite{BV-Paper2}\,. Here, we have written everything referring to the concrete case of the spectral fractional Laplacian in order to keep the exposition clear and to focus on the main ideas, but the arguments are devised in view of the wider applicability.

\medskip

\noindent {\sc Existence theory.}  The paper is complemented with a brief presentation of the theory of existence and uniqueness of the class of very weak solutions that we use, see Section \ref{sect.Exist.Uniq}.  We refer to \cite{DPQRV1,DPQRV2} for the basic theory on existence, uniqueness and $L^q$ and $C^\alpha$ regularity. We refer to \cite{BV2012} for quantitative estimates for very weak solutions as well as the initial trace problem. We refer to \cite{BV2012,DPQRV1,DPQRV2} for a physical motivation and relevance of this nonlocal model.

We refer the reader to the final section of the paper for  comments, consequences, ideas on extensions, and related works.

\section{Preliminaries}\label{sect.prelim}
Before addressing the derivation of the quantitative estimates, which is the main goal of the paper, we review the basic facts needed to set up the problem. First, we review some notations. The symbol $\infty$ will always denote $+\infty$. We will write $a\wedge b=\min\{a,b\}$. We will use the notation $f\asymp g$ if and only if there exists constants $c_0,c_1>0$ such that $c_0\,g\le f\le c_1 g$\,. We denote by $\LL^p_{\Phi_1}(\Omega)$ the weighted $\LL^p$ space $\LL^p(\Omega\,,\, \Phi_1\dx)$, endowed with the norm
\[
\|f\|_{\LL^p_{\Phi_1}(\Omega)}=\left(\int_{\Omega} |f(x)|^p\Phi_1(x)\dx\right)^{\frac{1}{p}}\,.
\]
We will always consider smooth domains $\Omega$ - with boundary at least $C^{1,1}$ - unless explicitly stated. Throughout the whole paper, we will consider $m>1$ and $0<s\le 1$\,. The a priori bounds include some exponents that will have the following values:
\[
\vartheta_1=\vartheta_{1,0}=\frac{1}{2s+N(m-1)}, \qquad\mbox{and}\qquad
\vartheta_{1,\gamma}=\frac{1}{2s+(N+\gamma)(m-1)}\,,
\]
where $\gamma>0$. By universal constant we mean a constant $K>0$ that may depend only on $N, m, s$ and $\Omega$\,, but not on $u$\,.

\subsection{Fractional Laplacian operators on bounded domains}\label{ssect.Def.Fract.Lapl}

Next, we discuss the various possible definitions of the concept of fractional Laplacian operator. When we consider operators and equations on the whole Euclidean space $\RR^N$, there is a natural concept of fractional Laplacian that can be defined in several equivalent ways: on the one hand, the definition of the nonlocal operator $(-\Delta_{\RR^N})^{s}$, known as the Laplacian of order $s$, is done by means of Fourier series
\begin{equation}\label{sLapl.Rd.Fourier}
  ((-\Delta_{\RR^N})^{s}f)^{\widehat{}}(\xi)=|\xi|^{2s} \hat f (\xi)\,,
\end{equation}
and can be used for positive and negative values of $s$, cf.  Stein, \cite{Stein70}. If $0<s<1$, we can also use the representation by means of an
hypersingular kernel,
\begin{equation}\label{sLapl.Rd.Kernel}
(-\Delta_{\RR^N})^{s}  g(x)= c_{d,s}\mbox{
P.V.}\int_{\mathbb{R}^N} \frac{g(x)-g(z)}{|x-z|^{d+2s}}\,dz,
\end{equation}
where $c_{N,s}>0$ is a normalization constant. Another classical way of defining the fractional powers of a linear self-adjoint nonnegative operator, is expressed in terms of the  semigroup associated to the standard Laplacian operator, which in our case reads
\begin{equation}\label{sLapl.Rd.Semigroup}
\displaystyle(-\Delta_{\RR^N})^{s}
g(x)=\frac1{\Gamma(-s)}\int_0^\infty
\left(e^{t\Delta_{\RR^N}}g(x)-g(x)\right)\frac{dt}{t^{1+s}}.
\end{equation}
All the above definitions are equivalent when dealing with the Laplacian on the whole space $\RR^N$.

When we consider the equation \eqref{FPME.equation} posed on bounded domains, the definition via the Fourier transform does not apply and different choices  appear as possible natural definitions of the fractional Laplacian.

\noindent $\bullet$ On one hand, starting  from the classical Dirichlet Laplacian $\Delta_{\Omega}$ on the domain $\Omega$\,, then the so-called {\em spectral definition} of the fractional power of $\Delta_{\Omega}$ uses the previous formula in terms of the semigroup associated to the Laplacian, namely
\begin{equation}\label{sLapl.Omega.Spectral}
\displaystyle(-\Delta_{\Omega})^{s}
g(x)=\sum_{j=1}^{\infty}\lambda_j^s\, \hat{g}_j\, \phi_j(x)
=\frac1{\Gamma(-s)}\int_0^\infty
\left(e^{t\Delta_{\Omega}}g(x)-g(x)\right)\frac{dt}{t^{1+s}}.
\end{equation}
where $\lambda_j>0$, $j=1,2,\ldots$ are the eigenvalues of the Dirichlet Laplacian on $\Omega$\,, written in increasing order and repeated according to their multiplicity, and $\phi_j$ are the corresponding normalized eigenfunctions and
\[
\hat{g}_j=\int_\Omega g(x)\phi_j(x)\dx\,,\qquad\mbox{with}\qquad \|\phi_j\|_{\LL^2(\Omega)}=1\,.
\]
 We will denote the operator defined in such a way as $\A_1=(-\Delta_{\Omega})^s$\,, and call it the \textit{spectral fractional Laplacian}, SFL for short. In this case, the initial and boundary conditions associated to the fractional diffusion equation \eqref{FPME.equation} read
\begin{equation}\label{FPME.Dirichlet.conditions.Spectral}
\left\{
\begin{array}{lll}
u(t,x)=0\,,\; &\mbox{in }(0,\infty)\times\partial\Omega\,,\\
u(0,\cdot)=u_0\,,\; &\mbox{in }\Omega\,,
\end{array}
\right.
\end{equation}
We refer to \cite{DPQRV1,DPQRV2} for the basic theory of weak solutions to the Dirichlet Problem \eqref{FPME.equation}--\eqref{FPME.Dirichlet.conditions.Spectral} with $u_0\in\LL^1(\Omega)$.

Let us list some properties of the operator:  $\A_1=(-\Delta_{\Omega})^s$ is a self-adjoint operator on $\LL^2(\Omega)$\,, with a discrete spectrum: we will denote by $\lambda_j^s>0$, $j=1,2,\ldots$ its eigenvalues (which are the family of $s$-power of the eigenvalues of the Dirichlet Laplacian), written in increasing order and repeated according to their multiplicity. We will denote by \ $\phi_{j}$ \ the corresponding normalized eigenfunctions which are exactly the same as the Dirichlet Laplacian, therefore they are as smooth as the boundary allows, namely when $\partial\Omega$ is $C^k$, then  $\phi_j\in C^{\infty}(\Omega)\cap C^k(\overline{\Omega})$\,.

The definition of the Fractional Laplacian via the Caffarelli-Silvestre extension \cite{Caffarelli-Silvestre} has been extended to the case of bounded domains by Cabr\'e and Tan \cite{Cabre-Tan}  by using as extended domain the cylinder ${\cal C}=(0,\infty)\times \Omega$ in $\RR^{N+1}_+$, and by putting zero boundary conditions on the lateral boundary of that cylinder. It is proved that this definition is equivalent to the SFL. See also \cite{capella-d-d-s}. \normalcolor

\noindent $\bullet$ On the other hand, we can define a fractional Laplacian operator by using the integral representation \eqref{sLapl.Rd.Kernel} in terms of hypersingular kernels and ``restrict'' the operator to functions that are zero outside $\Omega$: we will denote the operator defined in such a way as $\A_2=(-\Delta_{|\Omega})^s$\,, and call it the \textit{restricted fractional Laplacian}, RFL for short. In this case, the initial and boundary conditions associated to the fractional diffusion equation \eqref{FPME.equation} read
\begin{equation}\label{FPME.Dirichlet.conditions.Restricted}
\left\{
\begin{array}{lll}
u(t,x)=0\,,\; &\mbox{in }(0,\infty)\times\RR^N\setminus \Omega\,,\\
u(0,\cdot)=u_0\,,\; &\mbox{in }\Omega\,,
\end{array}
\right.
\end{equation}
We refer to \cite{SV1} for a more detailed discussion (and references) about the differences between the Spectral and the Restricted fractional Laplacian. The authors of \cite{SV1} call the second type simply Fractional Laplacian, but we feel that the absence of descriptive name leads to confusion.


\subsection{Functional setup}\label{ssect.Funct.Setup}

We will devote this paper to study Problem \ref{FPME.equation}--\ref{FPME.Dirichlet.conditions.Spectral} with the fractional Laplacian operator $\A_1$ that we write in the sequel simply $\A$\,. We denote by $\lambda_k$  the  eigenvalues of $\A$, written in increasing order and repeated according to their multiplicity, and the $\Phi_k$ being the corresponding normalized eigenfunctions respectively. It is  known that the first eigenfunction $\Phi_1$ is positive. The boundary behaviour of $\Phi_1$ is as follows: when the domain is $C^{1,1}$ smooth the following estimate holds
\begin{equation}\label{Phi1.est}
\Phi_{1}(x)
\asymp \big(\dist(x, \partial\Omega)\wedge 1\big)\qquad\mbox{for all }x\in\Omega\,,
\end{equation}
and  moreover $\Phi_1$ is $C^1$ up to the boundary. It is well known that $\A:H(\Omega)\to H^*(\Omega)$ is an isomorphism between the Hilbert space $H(\Omega)$ and its dual $H^*(\Omega)$, where
\begin{equation}\label{def.H.space}
H:=\left\{f\in \LL^2(\Omega)\,\Big|\, \sum_{k=1}^{+\infty}\lambda_k \hat{f}_k^2<+\infty\right\}\,, \quad\mbox{where}\quad \hat{f}_k=\int_\Omega f(x)\Phi_k(x)\dx\,.
\end{equation}
We will use the following norms on $H$ and $H^*$ respectively. Firstly,
\begin{equation}\label{H.norms}
\|f\|_{H(\Omega)}^2=\sum_{k=1}^{+\infty}\lambda_k \hat{f}_k^2= \|\AM f\|_{\LL^2(\Omega)}^2 = \int_\Omega f\A f\dx\,,
\end{equation}
and, letting $\AI: H^*\to H$ be the inverse of the isomorphism $\A$, we have
\begin{equation}\label{H*.norms}
\|f\|_{H^*(\Omega)}^2=\sum_{k=1}^{+\infty}\lambda_k^{-1} \hat{f}_k^2= \|\AIM f\|_{\LL^2(\Omega)}^2 = \int_\Omega f\AI f\dx\,.
\end{equation}
Notice that $H^*(\Omega)$ contains all the functions $f\in \LL^1_{\Phi_1}(\Omega)\cap\LL^\infty(\Omega)$: here is a sketch of the proof for $f\ge 0$
\[
\int_\Omega f\AI f\dx \le \left\|\frac{f}{\Phi_1}\right\|_{\LL^\infty}\int_\Omega \Phi_1  \AI f\dx
=\left\|\frac{f}{\Phi_1}\right\|_{\LL^\infty}\int_\Omega f\AI \Phi_1\dx=\lambda_1^{-1}\left\|\frac{f}{\Phi_1}\right\|_{\LL^\infty} \int_\Omega f\Phi_1\dx\,.
\]
For our choice of fractional Laplacian operator, it is also known that  we can identify the space $H(\Omega)$ with more usual fractional Sobolev spaces:
\begin{equation}\label{H-Hs}
H(\Omega)=\left\{\begin{array}{lll}
H^s_0(\Omega)\,, &\qquad\mbox{if }\frac{1}{2}<s\le 1\,,\\[2mm]
H^{1/2}_{00}(\Omega)\,,&\qquad\mbox{if }s=\frac{1}{2}\,,\\[2mm]
H^s(\Omega)\,, &\qquad\mbox{if }0<s< \frac{1}{2}\,.
\end{array}\right.
\end{equation}
Therefore the dual $H^*(\Omega)=H^{-s}(\Omega)$ if $s\ne 1/2$\,, and $H^{-1/2}(\Omega)\subset (H^{1/2}_{00}(\Omega))^*$ for $s=1/2$\,. Recall that $H^s_0$ is the closure of $C_c^\infty(\Omega)$ in the $H^s$-norm, which is equivalent to the $H$-norm above. Moreover, for $0< s\le 1/2$ we have the identity $H^s(\Omega)=H^s_0(\Omega)$\,, while both spaces are different for $s>1/2$\,. Finally, the following characterization holds:
\begin{equation}\label{def.H1200.2}
H^{\frac{1}{2}}_{00}(\Omega):= \left\{f\in H^{\frac{1}{2}}_0(\Omega)\;\big|\; \widetilde{\rm \delta}^{-\frac{1}{2}}f\in \LL^2(\Omega)\right\}\,,
\end{equation}
where $\widetilde{\rm \delta}$ is a smooth extension of the distance to the boundary\,. Moreover, we can equip $H^{\frac{1}{2}}_{00}(\Omega)$ with the following norm
\begin{equation}\label{norm.H1200.1}
\|f\|_{H^{\frac{1}{2}}_{00}(\Omega)}:= \left(\|f\|^2_{H^{\frac{1}{2}}(\Omega)}
    +\left\|{\widetilde{\delta}}^{-\frac{1}{2}}f\right\|_{\LL^2(\Omega)}^2\right)^{\frac{1}{2}}\,.
\end{equation}
Finally, the space $H^s_{00}(\Omega)$ is strictly included in $H^s_{0}(\Omega)$ with a strictly finer topology. The dual space of $H_{00}^{\frac{1}{2}}(\Omega)\subsetneqq  H_0^{\frac{1}{2}}(\Omega)$ satisfies
$$
H^{-\frac{1}{2}}(\Omega)= \left(H^{-\frac{1}{2}}_0(\Omega)\right)^*= \left(H^{-\frac{1}{2}}(\Omega)\right)^*
\subsetneqq \left(H_{00}^{\frac{1}{2}}(\Omega)\right)^*\,.
$$
For more details on the above discussion we refer to \cite{Adams2003, LM}.

\medskip

\noindent\textbf{The inverse of $\A$ and Green functions estimates. }The Green operator $\AI:H^*(\Omega)\to H(\Omega)$ is defined as the inverse of the operator $\A$, and it can be shown that it has a symmetric kernel $G_\Omega(x,y)$, which is the Green function:
\begin{equation}\label{def.Green.0}
\AI f(x_0):=\sum_{k=1}^{+\infty}\lambda_k^{-1} \hat{f}_k\Phi_k(x_0)=\int_\Omega G_\Omega(x,x_0)f(x)\dx\,.\normalcolor
\end{equation}
We provide further details in Appendix \ref{SSection.Green}\,. As already mentioned in the introduction, one of the novelties of this paper is represented by the technique used in the proof of the lower and upper estimates. The basic and main ingredient is the knowledge of good estimates for the Green function.\normalcolor

We now state the two kinds of estimates for the Green function of $\A$ that will be used in the proofs of our main results. They are Type I: for all $x,x_0\in \Omega$ we have
\begin{equation}\label{typeI.Green.est}
 0\le G_\Omega(x,x_0)\le \frac{c_{1,\Omega}}{|x-x_0|^{N-2s}}\,,
\end{equation}
and Type II:
\begin{equation}\label{typeII.Green.est}
 c_{0,\Omega}\Phi_1(x)\Phi_1(x_0)\le G_{\Omega}(x,x_0)\le
\frac{c_{1,\Omega}}{|x-x_0|^{N-2s}}
\left(\frac{\Phi_1(x)}{|x-x_0| }\wedge 1\right)
\left(\frac{\Phi_1(x_0)}{|x-x_0| }\wedge 1\right)\,.
\end{equation}
Recall that $\Phi_1$ is the first eigenfunction of $\A$\,, that satisfies estimates \eqref{Phi1.est}; recall also that $(a\wedge 1)(b\wedge 1)=\min\{1,a,b,ab\}$ for all $a,b\ge 0$\,.

A key point of the Type I estimate \eqref{typeI.Green.est} is the integrability property of the Green function. Indeed, by the H\"older inequality, the estimate $0\le G_{\Omega}(x,x_0)\le c_{1,\Omega}/|x-x_0|^{N-2s}$ implies that $G_{\Omega}(\cdot,x_0)\in \LL^r(\Omega)$ with $r<d/(N-2s)$\,, so that
\begin{equation}\label{Intro.1}
\int_{\Omega} \varphi\,G_{\Omega}(x,x_0)\dx \le \|\varphi\|_{r'}\|G_{\Omega}(\cdot,x_0)\|_r\,,\qquad\mbox{with }r'=\frac{r}{r-1}
\end{equation}
and clearly $r'>N/(2s)$. In our proofs, we will need more precise estimates both from above and from below, which also take into account the boundary behaviour (through Type II estimates), and this will be carefully done in Lemmata \ref{Lem.Green} and \ref{Lem.Green.2} of Appendix \ref{SSection.Green}\,.

Type I and Type II estimates are essentially used  in the sequel to obtain our main lower and upper estimates. On one hand,  the Type I estimate gives global absolute bounds from above and $\LL^1-\LL^\infty$ smoothing effects. On the other hand, the Type II estimate will imply sharp lower and upper boundary behaviour. Notice that the lower bound of Type II \eqref{typeII.Green.est} is weaker than the best known estimate on Green functions \eqref{Gree.est.0}, but such weaker lower bounds have the advantage that may hold in more general settings, see \cite{D1,D2,Davies1,Davies2,DS}\,, which may be useful for further research.


\section{Precise statement of the problem}\label{sect.statement}

We will consider the homogeneous Dirichlet problem
\begin{equation}\label{FPME.prob}
\left\{\begin{array}{lll}
u_t+\A(u^m)=0 &  ~ {\rm in}~ (0,+\infty)\times \Omega\\
u(0,x)=u_0(x) & ~{\rm in}~ \Omega \\
u(t,x)=0 & ~\mbox{ on $(0,+\infty)\times\partial \Omega$.}
\end{array}
\right.
\end{equation}
where $\Omega\subset \RR^N$ is a bounded domain with smooth boundary, $m>1$, $0<s<1$ and $N\ge 1$\,; we recall that the linear operator $\A$ will be in principle the spectral fractional Laplacian $\A_1=(-\Delta_{\Omega})^{s}$ defined above. The boundary conditions are interpreted as the standard Dirichlet conditions for the extended Caffarelli-Silvestre function $W(x,y,t)$, which is defined in the cylinder $\mathcal{C}=\Omega\times \RR^+$\,, cf. \cite{Cabre-Tan, DPQRV1,DPQRV2}\,.

Let us recall the definitions of weak and strong solution used in \cite{DPQRV1,DPQRV2}:\normalcolor
\begin{defn}\label{Def.Weak.Sol} A function $u$ is a {\sl weak} solution to Equation \eqref{FPME.equation} if:
\begin{itemize}
\item $u\in C((0,\infty): \LL^1(\Omega))$\,, $|u|^{m-1}u \in \LL^2_{\rm loc}\left((0,\infty):H(\Omega)\right)$;
\item  The identity
\begin{equation}
\displaystyle \int_0^\infty\int_{\Omega}u\dfrac{\partial \psi}{\partial t}\,\dx\dt
-\int_0^\infty\int_{\Omega} \left(\AM|u|^{m-1}u\right)\,\left(\AM\psi\right)\,\dx \dt=0.
\end{equation}
holds for every $\psi\in C^1_0((0,+\infty)\times \Omega)$\,.\normalcolor
\item A {\sl weak} solution to the Cauchy-Dirichlet problem \eqref{FPME.prob} is a weak solution to Equation \eqref{FPME.equation} such that moreover $u\in C([0,\infty): L^1_{\Phi_1}(\Omega))$ and $u(\cdot,0)=u_0$ almost everywhere in $\Omega$.
\item A weak solution to the Cauchy-Dirichlet problem \eqref{FPME.prob} is a {\sl strong} solution if moreover\\ $\partial_t u$ and $\A u^m \in~\LL^\infty((\tau,\infty): L^1(\Omega))$ for every $\tau>0$.
\end{itemize}
\end{defn}
Existence and uniqueness of strong solutions has been proved in \cite{DPQRV2} as follows:
\begin{thm}\label{th:general}
Let $m>1$ and $s\in (0,1]$. For every $u_0\in \LL^1(\Omega)$ there exists a unique strong solution of
the Dirichlet Problem $\eqref{FPME.prob}$. The strong solutions enjoy a number of properties, like the maximum principle. Nonnegative solutions are actually continuous and positive everywhere. Moreover, $u\in \LL^\infty([0,\infty):\LL^1(\Omega))$ and also $u\in \LL^\infty((\tau,\infty)\times\LL^\infty(\Omega))$ for all $\tau>0$\,.
\end{thm}

We can enlarge the concept of solution in the following way
\begin{defn}\label{Def.Very.Weak.Sol} A function $u$ is a {\sl very weak} solution to Equation \eqref{FPME.equation} if:
\begin{itemize}
\item $u\in C((0,\infty): \LL^1_{\Phi_1}(\Omega))$\,, $|u|^{m-1}u \in \LL^1\left((0,\infty):L^1_{\Phi_1}(\Omega)\right)$;
\item  The identity
\begin{equation}
\displaystyle \int_0^\infty\int_{\Omega}u\dfrac{\partial \psi}{\partial t}\,\dx\dt
-\int_0^\infty\int_{\Omega} |u|^{m-1}u\,\A\psi\,\dx \dt=0.
\end{equation}
holds for every $\psi\in C^2_c((0,+\infty)\times\overline{\Omega})$\,.
\item A {\sl very weak} solution to the Cauchy-Dirichlet problem \eqref{FPME.prob} is a very weak solution to Equation \eqref{FPME.equation} such that moreover $u\in C([0,\infty): L^1_{\Phi_1}(\Omega))$.
\end{itemize}
\end{defn}
In the terminology that is used in the specialized literature this is called a {\em weighted very weak solution}. The use of the weighted space is natural in view of our calculations below. In the elliptic case, the above definition of very weak solution, has been first given by Brezis (unpublished note), mentioned and used in \cite{DR09}, see also \cite{BCMR96}\,.\normalcolor

We will use a rather abstract approach in the derivation of our estimates. With this in mind, it will be convenient to use a still more general notion of solution. Thus, we introduce the concept of \textit{weak dual solution}, expressed in terms of the  problem involving the inverse of the fractional Laplacian.

\begin{defn}\label{Def.Very.Weak.Sol.Dual} A function $u$ is a {\sl weak dual} solution to Equation \eqref{FPME.equation} in $Q_T=(0,T)\times \RR^N$ if:
\begin{itemize}
\item $u\in C((0,T): \LL^1_{\Phi_1}(\Omega))$\,, $|u|^{m-1}u \in \LL^1\left((0,T):\LL^1_{\Phi_1}(\Omega)\right)$;
\item  The identity
\begin{equation}
\displaystyle \int_0^T\int_{\Omega}\AI (u) \,\dfrac{\partial \psi}{\partial t}\,\dx\dt
-\int_0^\infty\int_{\Omega} |u|^{m-1}u\,\psi\,\dx \dt=0.
\end{equation}
holds for every test function $\psi$ such that  $\psi/\Phi_1\in C^1_c((0,T): \LL^\infty(\Omega))$\,.
\item A {\sl weak dual} solution to the Cauchy-Dirichlet problem \eqref{FPME.prob} is a  weak dual solution to Equation \eqref{FPME.equation} such that moreover $u\in C([0,T): \LL^1_{\Phi_1}(\Omega))$  and $u(0,x)=u_0\in \LL^1_{\Phi_1}(\Omega))$.
\end{itemize}
\end{defn}
\noindent\textbf{Remarks. } (i) Roughly speaking, we are considering the weak solution to the equation $\partial_t U= u^m$\,, where $U=\AI u$\,, posed on the bounded domain $\Omega$ with homogeneous Dirichlet conditions. This ``dual formulation'' is helpful, and it has been introduced in \cite{Vaz2012} in the case the Cauchy problem for the porous medium equation on the whole $\RR^N$\,, for any $0<s\le 1$, and previously by Pierre \cite{Pierre} in the local case $s=1$.\\
(ii) In the elliptic case, this type of definition of very weak solution has been given in \cite{KR07}, for $s=1$.\normalcolor\\
(iii) The finite existence time $T>0$ is used in the definition for generality, but the we use below for simplicity the choice $T=\infty$ without loss of generality in the arguments.\\
(iv) Notice that the condition $\psi/\Phi_1\in C^1_c((0,T): \LL^\infty(\Omega))$ implies that $\|\psi(t,\cdot)/\Phi_1\|_{\LL^\infty(\Omega)}<+\infty$  and $\|\partial_t\psi(t,\cdot)/\Phi_1\|_{\LL^\infty(\Omega)}<+\infty$ for all $t\ge 0$ and moreover, are compactly supported functions of $t>0$ therefore in $\LL^1(0,\infty)$.\\
(v) It can be shown that strong solutions are weak dual solutions. \\
(vi) Since the nonlinearity has a power form, it is easy to see that the set of very weak  solutions defined globally in time, $T=\infty$, has the property of scale invariance, i.e., if $u(t,x)$ is a solution, then so is the function $\widetilde u$ defined as \ $\widetilde u(t,x)=\lambda^{m-1} u(\lambda t, x)$ for every constant $\lambda>0$.

In our study we will need a somewhat more restricted class of solutions.
\begin{defn}\label{Def.Very.Weak.Sol.Dual.2} We consider a class $\mathcal{S}$ of nonnegative  weak dual solutions $u$ to the Dirichlet problem \eqref{FPME.prob} with initial data in $u_0\in\LL^1_{\Phi_1}(\Omega)$\,, such that  $(i)$ the map $u_0\mapsto u(t)$ is  order preserving in $ \LL^1_{\Phi_1}(\Omega)$; $(ii)$ the class is also scale invariant; $(iii)$ for all $t>0$ we have $u(t)\in L^p(\Omega)$ with $p>N/2s$.
\end{defn}

The results of sections \ref{sec.1set} to \ref{sect.Elliptic.Prob} hold in the greater generality of operators for which there is a class of solutions $\cal S$, for which the Green function representation \eqref{def.Green.0} holds as well as the Type I and II estimates that have been mentioned.

The reader may prefer to concentrate on the application to the SFL in what follows. We will devote Section \ref{sect.Exist.Uniq} below to the construction of solutions that build up a class $\mathcal{S}$ satisfying these conditions in the case of the Spectral Fractional Laplacian, SFL. The class contains all nonnegative strong solutions.

\begin{thm}\label{thm.L1weight.exist}
For every  nonnegative $ u_0\in\LL^1_{\Phi_1}(\Omega)$ there exists a unique minimal very weak solution to the Dirichlet problem \eqref{FPME.prob} with ${\cal L}=(-\Delta_{\Omega})^{s}$, the SFL. Such a solution is obtained as the monotone limit of the strong solutions that exist and are unique when the initial data are bounded. The minimal very weak solution is continuous in the weighted space $u\in C([0,\infty):\LL^1_{\Phi_1}(\Omega))$. The very weak solutions are weak dual solutions and the set of such solutions has the properties needed to form a class of type~$\mathcal{S}$.
\end{thm}

\section{First set of estimates}\label{sec.1set}

This section is devoted to obtaining some basic estimates which will be essential in the proof of the main results. In particular, the pointwise bounds of Proposition \eqref{thm.NLE.PME} will imply both lower, upper and boundary estimates, that will combine into Harnack inequalities.

The first property we examine is the almost  monotonicity in time.  It is based on the famous B\'enilan-Crandall differential estimate, which depends only on the homogeneity of the equations, scale invariance,  and the property of order preservation, cf. \cite{BCr}.

\begin{lem}\label{lem.BC.est} For every solution $u\ge 0$  in the class $\cal S$ we have the differential estimate
\begin{equation}
u_t\ge -\frac{u}{(m-1)t}
\end{equation}
in the sense of distributions in $Q$. We recall that $m>1$.
\end{lem}

In the case of the spectral Dirichlet Laplacian, the proof has been verified in \cite{DPQRV2}, and the argument applies here. This property means that for almost all $(x,t)$ the function
$u(x,t)t^{1/(m-1)}$ is non-decreasing in time for fixed $x$. Therefore, for all $\tau\le t_1$:
\begin{equation}\label{BC.est.PME}
 u(\tau,x)\le \left(\frac{t_1}{\tau}\right)^{\frac{1}{m-1}}u({t_1},x)\,, \mbox{ a.e. in }\Omega\,.
\end{equation}
Let us now establish some integral and pointwise estimates.  We denote by $G_\Omega$  the Green function of $\A$, as described in Section \ref{ssect.Funct.Setup}, see more in  Section \ref{SSection.Green}.

\begin{prop}\label{thm.NLE.PME}  let $u\ge 0$ be a solution in the class $\mathcal{S}$ of very weak solutions to the Dirichlet problem \eqref{FPME.prob}. Then,
\begin{equation}\label{thm.NLE.PME.estim.0}
\int_{\Omega}u(t,x)G_{\Omega}(x , x_0)\dx\le \int_{\Omega}u_0(x)G_{\Omega}(x , x_0)\dx \qquad\mbox{for all $t> 0$\,.}
\end{equation}
Moreover, for almost every $0< t_0\le t_1 \le t$ and almost every $x_0\in \Omega$\,, we have
\begin{equation}\label{thm.NLE.PME.estim}
\left(\frac{t_0}{t_1}\right)^{\frac{m}{m-1}}(t_1-t_0)\,u^m(t_0,x_0) \le \int_{\Omega}u(t_0,x)G_{\Omega}(x , x_0)\dx - \int_{\Omega}u({t_1},x)G_{\Omega}(x , x_0)\dx \le (m-1)\frac{t^{\frac{m}{m-1}}}{t_0^{\frac{1}{m-1}}} \,u^m(t,x_0)\,.
\end{equation}
\end{prop}

\noindent\textbf{Remarks. }(i) A very weak solution  $u(x,t)\ge 0$ belonging to the class $\mathcal{S}$ of solutions to the Dirichlet problem \eqref{FPME.prob} has the property that $u(t)\in \LL^p(\Omega)$ with $p>N/2s$\,, therefore $\int_{\Omega}u(t,x)G_{\Omega}(x , x_0)\dx<+\infty$ for all $t>0$\,, simply as a consequence of H\"older inequality, as already mentioned in \eqref{Intro.1}. Therefore, using the pointwise inequality \eqref{thm.NLE.PME.estim} we may conclude that $u(t)$ is bounded for all $t>0$. On the other hand, we only assume that  $0\le u_0\in \LL^1_{\Phi_1(\Omega)}$, hence $\int_{\Omega}u_0(x)G_{\Omega}(x , x_0)\dx$ may be unbounded.\\
(ii) The most relevant part of the last estimate for the applications  is represented by the middle term: we will use it to deduce the lower and upper estimates, that will then be combined into quantitative Harnack inequalities.

\medskip

\noindent {\sl Proof.~}We split the proof in several steps.

\noindent \textsc{Step 1. }\textit{$\LL^1$-weighted estimates. } We will use the definition \eqref{Def.Very.Weak.Sol.Dual} of very weak solution, with a test function of the form $\psi(t,x)=\psi_1(t)\psi_2(x)$, where $\psi_1(t)\in  C^1_c((0,+\infty))$ and $\psi_2/\Phi_1\in \LL^\infty(\Omega)$. It follows that $u\in C((0,\infty): \LL^1_{\Phi_1}(\Omega))$, $u^m \in \LL^1\left((0,\infty):\LL^1_{\Phi_1}(\Omega)\right)$ satisfies the identity
\begin{equation}\label{step.1.thm.repr.0}\begin{split}
\int_0^{+\infty}\psi'_1(\tau)\int_{\Omega}u(\tau,x)  \AI\psi_2(x)\dx\rd\tau
 &= \int_0^{+\infty}\psi_1(\tau)\int_{\Omega}u^m(\tau,x)\,\psi_2(x)\dx \rd\tau\,,
\end{split}
\end{equation}
where in the left-hand side we have used the symmetry of the operator $\AI$. Notice that the space integral on the left-hand side of the formula is bounded, because of the following argument: we write $\psi_2= v\Phi_1$ with $v(x)$ bounded, and recall that $u\ge 0$. Then,
\[
\int_\Omega \AI u(t)\, \psi_2\dx \le \left\|v\right\|_{\LL^\infty}\int_\Omega \Phi_1  \AI u(t)\dx
=\left\|v\right\|_{\LL^\infty}\int_\Omega u(t)\AI \Phi_1\dx=\lambda_1^{-1}\left\|v\right\|_{\LL^\infty} \int_\Omega u(t)\Phi_1\dx\,.
\]
We now want to pass to the limit in \eqref{step.1.thm.repr.0} and prove that for all $0\le t_0\le t_1$ and for all $\psi_2(x)$\,, with $\psi_2: \overline{\Omega}\to \RR$ measurable and $\|\psi_2/\Phi_1\|_{\LL^\infty(\Omega)}<+\infty$\,, we have
\begin{equation}\label{step.1.thm.repr}\begin{split}
\int_{\Omega}u(t_0,x)\AI \psi_2(x)\dx - \int_{\Omega}u({t_1},x)\AI \psi_2(x)\dx=\int_{t_0}^{t_1}\int_{\Omega}u^m(\tau,x)\psi_2(x)\dx \rd\tau\,.
\end{split}
\end{equation}
This is rather standard: we only need to take $\psi_1(\tau)=\chi_{[t_0,t_1]}(\tau)$ as test function in formula \eqref{step.1.thm.repr.0}\,, so that $\psi_1'(\tau)=\delta_{t_0}(\tau)-\delta_{t_1}(\tau)$; this can be jusfified by considering a smooth approximation $\psi_{1,n}\in C_c^{\infty}(0,+\infty)$ so that $\psi_{1,n}\to \chi_{[t_0,t_1]}(\tau)$ in $\LL^\infty(0,+\infty)$\,, and so that $\psi'_{1,n} \to \delta_{t_0}(\tau)-\delta_{t_1}(\tau)$ in the sense of Radon measures with compact support. Clearly, these approximations are admissible test functions such that $\psi_{n}/\Phi_1\in C^1_c((0,+\infty): \LL^\infty(\Omega))$\,. Under the above assumptions, it is clear that
\[
\int_0^{+\infty}\psi'_{1,n}(\tau)\int_{\Omega}u(\tau,x)  \AI\psi_2(x)\dx\rd\tau \xrightarrow[n\to\, \infty]{}\int_{\Omega}u(t_0,x)\AI \psi_2(x)\dx -\int_{\Omega}u({t_1},x)\AI \psi_2(x)\dx\,,
\]
since $u\in C((0,\infty): \LL^1_{\Phi_1}(\Omega))$ implies that $\int_{\Omega}u({t_1},x)\AI \psi_2(x)\dx \in C^0(0,\infty)$\,, and since $\psi'_{1,n} \to \delta_{t_0}(\tau)-\delta_{t_1}(\tau)$  in the sense of Radon measures with compact support. On the other hand,
\[
\int_0^{+\infty}\psi_{1,n}(\tau)\int_{\Omega}u^m(\tau,x)\,\psi_2(x)\dx \rd\tau \xrightarrow[n\to\, \infty]{}
\int_{t_0}^{t_1}\int_{\Omega}u^m(\tau,x)\psi_2(x)\dx \rd\tau
\]
since $|u|^{m-1}u \in \LL^1\left((0,\infty):\LL^1_{\Phi_1}(\Omega)\right)$ implies that $\int_{\Omega}u^m(\tau,x)\psi_2(x)\dx \in \LL^1(0,\infty)$\,, and $\psi_{1,n}\to \chi_{[t_0,t_1]}(\tau)$ in $\LL^\infty(0,+\infty)$\,.

\noindent {\sc Note.} The above approximation also justifies the formal time-derivation of the $\LL^1_{\Phi_1}$-norm used in the proof of the Theorems \ref{thm.L1weight.PME} and \ref{thm.L1weight.contr}\,.

\noindent$\bullet~$\noindent\textsc{Step 2. }\textit{Proof of \eqref{thm.NLE.PME.estim.0}. }From \eqref{step.1.thm.repr} we  prove estimate \eqref{thm.NLE.PME.estim.0}, by first fixing $x_0\in \Omega$ and then taking a sequence of nonnegative test functions $\psi_{2,n}^{(x_0)}$ with $\psi_{2,n}^{(x_0)}: \overline{\Omega}\to \RR$ measurable and $\|\psi_2/\Phi_1\|_{\LL^\infty(\Omega)}<+\infty$\,, such that $\psi_{2,n}^{(x_0)}\to \delta_{x_0}$ as $n\to \infty$\,, in the sense of Radon measures. Therefore, $\AI \psi_{2,n}^{(x_0)}\to G_{\Omega}(\cdot , x_0)$ in $\LL^q(\Omega)$ for all $0<q<N/(N-2s)$, so that, taking $p=q/(q-1)>N/(2s)$\,, we have that
\[
\left|\int_{\Omega}u(\tau,x)  \AI\psi_{2,n}(x)\dx-\int_{\Omega}u(\tau,x) G_{\Omega}(x, x_0)\dx\right|\le \|u(\tau)\|_{\LL^p(\Omega)}\|\AI\psi_{2,n}-G_{\Omega}(\cdot , x_0)\|_{\LL^q(\Omega)}\to 0
\]
as $n\to\infty$\,,  and we recall that  $u\in \mathcal{S}$\,, therefore $u(t)\in \LL^p(\Omega)$ for all $t> 0$\,, with $p>N/(2s)$\,. Since the right-hand side of \eqref{step.1.thm.repr} is non-negative, we have proved the first estimate of the Theorem, \eqref{thm.NLE.PME.estim.0}\,.

\noindent$\bullet~$\noindent\textsc{Step 3. } We will use the monotonicity property \eqref{BC.est.PME} to estimate the right-hand side of  identity \eqref{step.1.thm.repr} from below and from above. More precisely, we will prove the following

\noindent {\bf Claim:} \textsl{For all $m>1$\,, for almost every $0\le t_0\le t_1 \le t$\,, and for all $\psi_2$ as in Step 1, we have:}
\begin{equation}\label{step.2.1.thm.repr}
\begin{split}
\left(\frac{t_0}{t_1}\right)^{\frac{m}{m-1}}(t_1-t_0)\int_\Omega u^m(t_0,x)\psi_2(x)\dx
&\le \int_{t_0}^{t_1}\int_\Omega u^m(\tau,x_0) \psi_2(x)\dx \rd\tau \\
&\le \frac{m-1}{t_0^{\frac{1}{m-1}}} t^{\frac{m}{m-1}}\,\int_\Omega u^m(t,x)\psi_2(x)\dx\,.
\end{split}
\end{equation}
Consider a smooth sequence  $\psi_{1,n}\in C_c^{\infty}(0,+\infty)$, $0\le \psi_{1,n}\le 1$, such that $\psi_{1,n}\to \chi_{[t_0,t_1]}$ in $\LL^\infty(0,+\infty)$ and such that $\supp(\psi_{1,n})\subseteq [t_0-1/n, t_1+1/n]$ and $\psi_{1,n}\ge \chi_{[t_0,t_1]}$\,. There are two cases.

\medskip

\noindent\textit{Upper estimates. }Let $n$ be so big that $0\le t_0-1/n \le t_1+1/n \le t$\,, and recall that $u\ge 0$\,, so that
\[
\begin{split}
\int_{0}^{\infty}\psi_{1,n}(\tau) &\int_{\Omega}u^m(\tau,x)\psi_2(x)\dx \rd\tau
\le  \int_{0}^{\infty}\psi_{1,n}(\tau) \left(\frac{t}{\tau}\right)^{\frac{m}{m-1}}\rd\tau \, \int_{\Omega}u^m(t ,x)\psi_2(x)\dx \\
&\le \|\psi_{1,n}\|_{\LL^\infty(0,+\infty)} \int_{t_0-1/n}^{t_1+1/n}\left(\frac{t }{\tau}\right)^{\frac{m}{m-1}}\rd\tau \, \int_{\Omega}u^m(t,x)\psi_2(x)\dx \\
&=(m-1)t^{\frac{m}{m-1}} \|\psi_{1,n}\|_{\LL^\infty(0,+\infty)} \left[\left(\frac{1}{t_0-\frac{1}{n}}\right)^{\frac{1}{m-1}}-\left(\frac{1}{t_1+\frac{1}{n}}\right)^{\frac{1}{m-1}}\right]\, \int_{\Omega}u^m(t,x)\psi_2(x)\dx \\
&\le (m-1)\frac{\|\psi_{1,n}\|_{\infty} }{\left(t_0-\frac{1}{n}\right)^{\frac{1}{m-1}}} t^{\frac{m}{m-1}}\, \int_{\Omega}u^m(t,x)\psi_2(x)\dx \,,
\end{split}
\]
where we have used inequality \eqref{BC.est.PME}\,, in the form  $u(\tau,x)\le \left(t/\tau\right)^{\frac{1}{m-1}}u(t,x)$ for all $t\ge t_1 +\frac{1}{n}\ge\tau$\,, and the assumptions on $\psi_{1,n}$. Let  $n\to \infty$ to get
\[
\int_{t_0}^{t_1}\int_{\Omega}u^m(\tau,x)\psi_2(x)\dx \rd\tau \le \frac{m-1}{t_0^{\frac{1}{m-1}}} t^{\frac{m}{m-1}}\, \int_{\Omega}u^m(t,x)\psi_2(x)\dx \qquad\mbox{for all $t\ge t_1 \ge t_0$}\,,
\]
since  $u^m \in \LL^1\left((0,\infty):\LL^1_{\Phi_1}(\Omega)\right)$ and  $\|\psi_{1,n}\|_{\LL^\infty(0,+\infty)}\to \|\chi_{[t_0,t_1]}\|_{\LL^\infty(0,+\infty)}=1$\,.

\noindent\textit{Lower estimates. }Let $n$ be so big that $0\le t_0-1/n \le t_1+1/n $. Since  $u\ge 0$\,, se have
\[
\begin{split}
\int_{0}^{\infty}\psi_{1,n}(\tau)\int_\Omega u^m(\tau,x)\psi_2(x)\dx \rd\tau
&\ge t_0^{\frac{m}{m-1}}\,\int_{t_0}^{t_1}\frac{\rd\tau}{\tau^{\frac{m}{m-1}}}\,\int_\Omega u^m(t_0,x)\psi_2(x)\dx \\
&=  (m-1){t_0^{\frac{m}{m-1}}} \left[\frac{1}{t_0^{\frac{1}{m-1}}}-\frac{1}{t_1^{\frac{1}{m-1}}}\right]\,\int_\Omega u^m(t_0,x)\psi_2(x)\dx
\end{split}
\]
where we have used monotonicity in the form $u(t_0-1/n,x)\le \left(\tau/(t_0-1/n)\right)^{\frac{1}{m-1}}u(\tau,x)$ for all $\tau \ge t_0 -\frac{1}{n}\ge0$\,, together with the fact that $\psi_{1,n}\ge0$ and $\psi_{1,n}\ge \chi_{[t_0,t_1]}$\,. Since the function $f(t)=t^{-\alpha}$ is convex for  $\alpha=1/(m-1)>0$, we may use the inequality $f(t)-f(t_0)\le f'(t)(t-t_0)$ to obtain
\[
\int_{0}^{\infty}\psi_{1,n}(\tau)\int_\Omega u^m(\tau,x)\psi_2(x)\dx\rd\tau
\ge \left(\frac{t_0}{t_1}\right)^{\frac{m}{m-1}}(t_1-t_0)\int_\Omega u^m(t_0,x)\psi_2(x)\dx\,.
\]
Letting now $n\to \infty$ we obtain
\[
\int_{t_0}^{t_1}\int_\Omega u^m(\tau,x)\psi_2(x)\dx\rd\tau
\ge \left(\frac{t_0}{t_1}\right)^{\frac{m}{m-1}}(t_1-t_0)\int_\Omega u^m(t_0,x)\psi_2(x)\dx\,,
\]
which is justified as before. The claim is proved.

\medskip

Summing up the results of the first  steps, \textsl{for all $m>1$\,, for every $0\le t_0\le t_1 \le t$\,, and for all $\psi_2(x)$\,, with $\psi_2: \overline{\Omega}\to \RR$ measurable and $\psi_2/\Phi_1$ bounded,  we have}
\begin{equation}\label{step.2.2.thm.repr}\begin{split}
\left(\frac{t_0}{t_1}\right)^{\frac{m}{m-1}}(t_1-t_0)\int_\Omega u^m(t_0,x)\psi_2(x)\dx
&\le\int_{\Omega}u(t_0,x)\AI \psi_2(x)\dx - \int_{\Omega}u({t_1},x)\AI \psi_2(x)\dx \\
&\le \frac{m-1}{t_0^{\frac{1}{m-1}}} t^{\frac{m}{m-1}}\,\int_\Omega u^m(t,x)\psi_2(x)\dx\,.
\end{split}\end{equation}

\noindent$\bullet~$\noindent\textsc{Step 4. } We will now prove formula \eqref{thm.NLE.PME.estim} by approximating the Green function $G_\Omega(x_0,\cdot)$ by means of a sequence of admissible test functions $\psi_{2,n}^{(x_0)}$.   For fixed $x_0\in\Omega$ we consider a sequence of test functions $\psi_{2,n}^{(x_0)}$ with $\psi_{2,n}^{(x_0)}: \overline{\Omega}\to \RR$ measurable and such that $\psi_{2,n}^{(x_0)}/\Phi_1$ is bounded\,, such that $\psi_{2,n}^{(x_0)}\to \delta_{x_0}$ as $n\to \infty$\,, in the sense of Radon measures. More specifically, we can choose $\psi_{2,n}^{(x_0)}(x)= |B_{1/n}(x_0)|^{-1}\,\chi_{B_{1/n}(x_0)}(x)$. Therefore, $\AI \psi_{2,n}^{(x_0)}\to G_{\Omega}(\cdot , x_0)$ in $\LL^q(\Omega)$ for all $0<q<N/(N-2s)$.   Then, for any fixed $\tau\ge 0$ we have:
\begin{equation}\label{step.3.1.thm.repr}
\lim_{n\to\infty}  \int_{\Omega}u^m(\tau,x)\psi_{2,n}^{(x_0)}(x)\dx
= \lim_{n\to\infty} |B_{1/n}(x_0)|^{-1}  \int_{B_{1/n}(x_0)}u^m(\tau,x)\dx
=   u^m(\tau,x_0)
\end{equation}
if $x_0$ is a Lebesgue point of the function $x\mapsto u(\tau,x)$\,; $u^m(\tau,x_0)$ is the corresponding Lebesgue value. If we apply this limit process at the points $\tau=t_0$ and $\tau=t_1$ we get for almost every $x_0$
\begin{equation}\label{step.3.3.thm.repr}\begin{split}
0\le \left(\frac{t_0}{t_1}\right)^{\frac{m}{m-1}} (t_1-t_0)u(t_0,x_0)
&\le\lim_{n\to\infty}\int_{\Omega}u(t_0,x)\AI \psi_{2,n}^{(x_0)}(x)\dx - \int_{\Omega}u({t_1},x)\AI \psi_{2,n}^{(x_0)}(x)\dx \\
&\le \frac{m-1}{t_0^{\frac{1}{m-1}}} t_1^{\frac{m}{m-1}}\,u(t_1,x_0)<+\infty\,,
\end{split}\end{equation}
Finally,  since $u\in \mathcal{S}$\,, then $u(t)\in \LL^p(\Omega)$ for all $t> 0$\,, with $p>N/(2s)$, and we have already seen in Step 2 that we have
\[
\int_{\Omega}u(t,x)\AI \psi_{2,n}^{(x_0)}(x)\dx \xrightarrow[n\to\, \infty]{} \int_{\Omega}u(t,x)G_{\Omega}(x,x_0)\dx\,,
\]
and formula \eqref{thm.NLE.PME.estim} follows for $t=t_1$\,. For $t$ larger than $t_1$ we use monotonicity. \qed


\section{Quantitative upper bounds}\label{Sect.Upper}

We use the pointwise lower  estimates of Theorem \ref{thm.NLE.PME} to prove two kinds of upper estimates: the smoothing effect and the absolute upper bound. The absolute bound holds up to the boundary and provides sharp estimates in terms of $\Phi_1\asymp \dist(\cdot,\partial\Omega)$\,.   The smoothing effect is sharp for small times, while the absolute upper bound is sharp for large times, see the remark at the end of this section for a more detailed explanation. Let  $m>1$ and  $u\ge 0$ is a solution in the class $\mathcal{S}$ of very weak solutions to the Dirichlet problem \eqref{FPME.prob}.

\subsection{Absolute bounds}

\begin{thm}[Absolute upper estimate]\label{thm.Upper.PME}
Let  $u$ be a solution in the class $\mathcal{S}$.  Then, there exists a universal constant $K_1>0$ such that the following estimates hold true
\begin{equation}\label{thm.Upper.PME.Absolute}
 \|u(t)\|_{\LL^\infty(\Omega)}\le\, K_1\,t^{-\frac{1}{m-1}}\,,\qquad\qquad\mbox{for all } \ t> 0\,.
\end{equation}
\end{thm}
This absolute bound proves a strong regularization which is independent of the initial datum. It also implies a sharp absolute boundary behaviour, as we shall see in the next theorem. On the other hand, the above absolute bounds \eqref{thm.Upper.PME.Absolute} hold for any $t> 0$\,, but they are not sharp for small times; precise bounds for small times are the smoothing effects of Theorem \ref{thm.Upper.1.PME}, where the upper bound depends on the initial datum.
The constant $K_1>0$ depends only on $N, m, s$ and $\Omega$\,, but not on $u$\,, and has an explicit form given in the proof.

\medskip

\noindent {\sl Proof.~} The proof is a consequence of the Type I bounds \eqref{typeI.Green.est} for the Green function.

\noindent$\bullet~$\textsc{Step 1. }\textit{Fundamental upper estimates. }We first recall the lower pointwise estimate of Theorem \ref{thm.NLE.PME}, that holds for any solution $u\in \mathcal{S}$:  for all $0\le t_0\le t_1 $ and $x_0\in \Omega$\,, we have that
\begin{equation}\label{Upper.PME.Step.1.1}
\left(\frac{t_0}{t_1}\right)^{\frac{m}{m-1}}(t_1-t_0)\,u^m(t_0,x_0) \le \int_{\Omega}u(t_0,x)G_{\Omega}(x , x_0)\dx - \int_{\Omega}u({t_1},x)G_{\Omega}(x , x_0)\dx\,.
\end{equation}
We choose $t_1=2t_0$ and recall that $u\ge0$\,, so that the above inequality \eqref{Upper.PME.Step.1.1} implies that
\begin{equation}\label{Upper.PME.Step.1.2}
u^m(t_0,x_0) \le \frac{2^{\frac{m}{m-1}}}{t_0}\int_{\Omega}u(t_0,x)G_{\Omega}(x , x_0)\dx\qquad\mbox{for all $t_0> 0$ and $x_0\in \Omega$}\,.
\end{equation}
This is a fundamental upper bound which encodes both the smoothing effect and the absolute upper bound, and therefore it is sharp both for large and for small times. A remarkable aspect of this upper bound is that it compares the $\LL^\infty$ norm and some integral norms at the same time $t_0> 0$\,.

\noindent The fact that $u\in \mathcal{S}$ guarantees that  $u(t)\in \LL^p(\Omega)$ for all $t> 0$\,, for some $p>N/(2s)$\,, so that
\[
\int_{\Omega}u(t_0,x)G_{\Omega}(x , x_0)\dx \le \|u(t_0)\|_{\LL^p(\Omega)}\,\|G_{\Omega}(\cdot , x_0)\|_{\LL^q(\Omega)}\le c_{2,\Omega}(q)\|u(t_0)\|_{\LL^p(\Omega)}<+\infty
\]
because $G_{\Omega}(\cdot , x_0)\in\LL^q(\Omega) $ for all $0< q<N/(N-2s)$\,, see Lemma \eqref{Lem.Green}. Therefore, we have
\begin{equation}\label{Upper.PME.Step.1.3}
u^m(t_0,x_0) \le c_{2,\Omega}(q)\frac{2^{\frac{m}{m-1}}}{t_0}\|u(t_0)\|_{\LL^p(\Omega)}\qquad\mbox{for all $t_0> 0$ and $x_0\in \Omega$}\,.
\end{equation}
so that $u(t_0)\in \LL^\infty(\Omega)$ for all $t_0>0$.

\noindent$\bullet~$\textsc{Step 2. }Let us estimate the right-hand side of the fundamental upper bound \eqref{Upper.PME.Step.1.2} in another way as follows:
\[
u^m(t_0,x_0) \le \frac{2^{\frac{m}{m-1}}}{t_0}\int_{\Omega}u(t_0,x)G_{\Omega}(x , x_0)\dx
\le \|u(t_0)\|_{\LL^\infty(\Omega)}\frac{2^{\frac{m}{m-1}}}{t_0}\int_{\Omega}G_{\Omega}(x , x_0)\dx\,.
\]
Therefore, taking the supremum over $x_0\in\Omega$ of both sides, we obtain:
\begin{equation}\label{Upper.PME.Step.2.1}
\|u(t_0)\|_{\LL^\infty(\Omega)}^{m-1}\le \frac{2^{\frac{m}{m-1}}}{t_0}\sup_{x_0\in\Omega}\int_{\Omega}G_{\Omega}(x , x_0)\dx
\le \frac{2^{\frac{m}{m-1}}c_{2,\Omega}}{t_0}:= \frac{K_1^{m-1}}{t_0}
\end{equation}
where we have used the bound $\sup\limits_{x_0\in\Omega}\int_{\Omega}G_{\Omega}(x , x_0)\dx \le c_{2,\Omega}$\,, given in Lemma \ref{Lem.Green}, inequality \eqref{Lem.Green.est.Upper.I}\,. Moreover, $K_1$ has the form
\begin{equation}\label{Upper.PME.Step.2.2}
K_1^{m-1}=2^{\frac{m}{m-1}}c_{2,\Omega}=2^{\frac{m}{m-1}}c_{1,\Omega}c_{N}\,\left(\diam(\Omega)+\frac{1}{s}\right)^N\|\Phi_1\|_{\LL^\infty(\Omega)}^{2  q}
\end{equation}
with $c_{1,\Omega}$ is the constant in the Type I estimates $\eqref{typeI.Green.est}$\,, and $c_N>0$ only depends on $N$\,. Note that $K_1$ blows up when $s\to 0^+$\,.\qed

\medskip

Next we  prove the sharp absolute upper boundary estimates. A key tool in the proof is the Integral Green function estimates of Lemma \ref{Lem.Green.2}. \normalcolor

\begin{thm}[Absolute boundary estimate]\label{thm.Upper.PME.II} There exists a universal constant $K_2>0$ such that the following estimate holds true
\begin{equation}\label{thm.Upper.PME.Boundary}
 u(t,x) \le K_2\, \frac{\Phi_1(x_0)^{\frac{1}{m}}}{t^{\frac{1}{m-1}}}\qquad\qquad\mbox{for  all $t> 0$ and $x\in \Omega$\,.}
\end{equation}
\end{thm}
Recall that $\Phi_1$ is the first eigenfunction of $\A$ and satisfies estimates \eqref{Phi1.est}\,.

\noindent {\sl Proof.~} It is a consequence of Type II bounds \eqref{typeII.Green.est} for the Green function, more precisely, we  use the Integral estimates II of Lemma \ref{Lem.Green.2}\,. We already know that $u(t)\in \LL^\infty(\Omega)$ for all $t>0$ by Theorem \ref{thm.Upper.PME}. Let us fix $t_0>0$, then the fundamental upper bound proved in Step 1 of Theorem \ref{thm.Upper.PME} reads: for almost all $x_0\in \Omega$
\begin{equation}\label{Upper.PME.Step.1.2.c}
u^m(t_0,x_0) \le \frac{2^{\frac{m}{m-1}}}{t_0}\int_{\Omega}u(t_0,x)G_{\Omega}(x , x_0)\dx
=\kappa_0\int_{\Omega}u(t_0,x)G_{\Omega}(x , x_0)\dx
\end{equation}
and this inequality guarantees that hypothesis \eqref{Lem.Green.2.hyp} of Lemma \ref{Lem.Green.2} of the Appendix holds true. Actually, it is proved that there exists a universal constant $c_{5,\Omega}>0$ such that
\begin{equation}
\mbox{if}\qquad u^m(x_0)\le \kappa_0\int_{\Omega} u(x)G_\Omega(x,x_0)\dx \qquad\mbox{then}\qquad
\int_{\Omega} u(x)G_\Omega(x,x_0)\dx\le c_{5,\Omega}^m\kappa_0^{\frac{1}{m-1}}\Phi_1(x_0)\,.
\end{equation}
Therefore,
\[
u^m(t_0, x_0)\le \kappa_0\int_{\Omega}u(t_0,x)G_{\Omega}(x , x_0)\dx\le c_{5,\Omega}^m \kappa_0^{\frac{m}{m-1}}\Phi_1(x_0)= c_{5,\Omega}^m 2^{\frac{m^2}{(m-1)^2}}\frac{\Phi_1(x_0)}{t_0^{\frac{m}{m-1}}}\,,
\]
which proves the desired upper bound \eqref{thm.Upper.PME.Boundary}, with universal constant $K_2=c_{5,\Omega} 2^{\frac{m}{(m-1)^2}}>0$\,,
where $c_{5,\Omega}>0$ is  given in Lemma \ref{Lem.Green.2}\,.\qed

\noindent\textbf{Remarks. } (i) The proof guarantees also the validity of the following inequality, which has its own interest and will be useful in the proof of the lower estimates of Section \ref{sect.positivity}:
\begin{equation}\label{thm.Upper.PME.Boundary.2}
\int_{\Omega}u(t_0,x)G_{\Omega}(x , x_0)\dx\le \overline{K}_2 \frac{\Phi_1(x_0)}{t_0^{\frac{1}{m-1}}}\,,
\qquad\qquad\mbox{for  all $t_0\ge 0$ and $x_0\in \Omega$\,,}
\end{equation}
where $\overline{K}_2=c_{5,\Omega}^m 2^{\frac{m}{(m-1)^2}}$ with $c_{5,\Omega}>0$ is as before\,.\\
(ii) The boundary behaviour is sharp. Indeed, we will obtain lower bounds with matching powers in Theorem \ref{thm.Lower.PME} of Section \ref{sect.positivity}. Also, the solutions obtained by separation of variables have the same boundary behaviour:
\[
U(t,x)=\frac{V(x)}{(t+h)^{\frac{1}{m-1}}}
\]
for any $h\in \RR$, where $V$ is a solution to the associated elliptic problem \eqref{Elliptic.prob} discussed in Section \ref{sect.Elliptic.Prob}, for which the sharp estimates of Theorem \ref{Thm.Elliptic.Harnack} apply: namely,  $V^m(x)\asymp \Phi_1(x)\asymp \big(\dist(x, \partial\Omega)\wedge 1\big)$ for all $x\in \Omega$\,.\\
 (iii) Consider the particular separation of variables solution $U(t,x)= V(x)t^{-\frac{1}{m-1}}$\,, which corresponds to the initial datum $U(0,x)=+\infty$\,. Such function has been called ``friedly giant'', see \cite{DK0,JLVmonats} since it represents the absolute maximum of the set of all nonnegative solutions at hand, i.e. solutions in the class $\mathcal{S}$. To be more precise, one can repeat the proof of Proposition 1.3 of \cite{JLVmonats} (since it basically depends  only on the validity of the above absolute bounds) to obtain that for all solution $u\in\mathcal{S}$ we have $u(t,x)\le U(t,x)$ for all $t>0$ and $x\in\Omega$\,. This has important consequences for the study of the asymptotic behaviour, see for example \cite{JLVmonats} for $s=1$, or \cite{BSV2013} for $s\in (0,1]$\,.


\subsection{The smoothing effects: $\LL^1$-$\LL^\infty$ and $\LL^1_{\Phi_1}$-$\LL^\infty$ }
Let $m>1$ and let $u\ge 0$ be a solution in the class $\mathcal{S}$ of very weak solutions to the Dirichlet problem \eqref{FPME.prob}, corresponding to the initial datum $0\le u_0\in \LL^1_{\Phi_1}(\Omega)$\,. We recall the exponents:
\[
\vartheta_{1,1}=\frac{1}{2s+(N+1)(m-1)}\qquad\mbox{and}\qquad \vartheta_1=\vartheta_{1,0}=\frac{1}{2s+N(m-1)}
\]
\begin{thm}\label{thm.Upper.1.PME}
 There exist universal constants $K_3,K_4>0$ such that the following estimates hold true.\\
\noindent\textsc{$\LL^1$-$\LL^\infty$ smoothing effect: }are consequence of Type I bounds \eqref{typeI.Green.est} for the Green function
\begin{equation}\label{thm.Upper.PME.Smoothing.1}
\|u(t)\|_{\LL^\infty(\Omega)}\le \frac{K_3}{t^{N\vartheta_1}}\|u(t)\|_{\LL^1(\Omega)}^{2s\vartheta_1}
 \le \frac{K_3 }{t^{N\vartheta_1}}\|u_0\|_{\LL^1(\Omega)}^{2s\vartheta_1} \qquad\qquad\mbox{for all $t> 0$\,.}
\end{equation}
\noindent\textsc{Intrinsic smoothing effect: }are consequence of Type II bounds \eqref{typeII.Green.est} for the Green function, for all $t> 0$ we have
\begin{equation}\label{thm.Upper.PME.Smoothing.2}
\|u(t)\|_{\LL^\infty(\Omega)}
    \le \frac{K_4}{t^{(N+1)\vartheta_{1,1}}}\|u(t)\|_{\LL^1_{\Phi_1}(\Omega)}^{2s\vartheta_{1,1}}
    \le \frac{K_4 }{t^{(N+1)\vartheta_{1,1}}}\|u_0\|_{\LL^1_{\Phi_1}(\Omega)}^{2s\vartheta_{1,1}}\,.
\end{equation}
\end{thm}

\noindent\textbf{Remarks. }(i) We have obtained weighted smoothing effects which are new to our knowledge also when $s=1$. Moreover, they apply to a class of nonnegative initial data $\LL^1_{\Phi_1}(\Omega) $ which is strictly larger than $\LL^1(\Omega)$.\\
(ii) Another novelty is represented by the fact that the smoothing effect occurs at the same time; this is new also when $s=1$\,.\\
 (iii) The label ``intrinsic'' for this smoothing effect is taken from Davies and Simon \cite{DS} who use the terminology {\sl intrinsic ultracontractivity} for estimates depending on the $\LL^2_{\Phi_1^2}$ norm. Indeed, our $\LL^1_{\Phi_1}$-$\LL^\infty$ smoothing effects, easily imply $\LL^2_{\Phi_1^2}$-$\LL^\infty$ smoothing effects, using the H\"older inequality.

\medskip

\noindent {\sl Proof.~}We already know from Theorem \ref{thm.Upper.PME} that the solutions under consideration are bounded for all $t>0$, here we are interested in understanding the behaviour for small times, where it enters in a substantial way the dependence on the initial datum. We split two cases, but we will use in both cases the following facts. Let $x_0\in \Omega$ and consider $B_r(x_0)$ with $r>0$ to be fixed later. Define the set $\Omega_r=\Omega\setminus\left(B_r(x_0)\cap \Omega\right)$ so that $\Omega\subseteq B_r(x_0)\cup \Omega_r$. Notice that the ball $B_r(x_0)$ need not to be included in $\Omega$\,. Then it is clear that $\forall x\in\Omega_r$ we have $|x-x_0|\ge r$\,.

\noindent$\bullet~$\textit{$\LL^1-\LL^\infty$ Smoothing effects via Type I estimates. } We will use the fundamental upper estimates \eqref{Upper.PME.Step.1.2} of Step 1 of the proof of Theorem \ref{thm.Upper.PME}\,, together with Type I estimates \eqref{typeI.Green.est}\,, namely $G_{\Omega}(x,x_0)\le c_{1,\Omega}|x-x_0|^{-(N-2s)}$\,, to obtain
\begin{equation*}\begin{split}
u^m (t_0,x_0) &\le \frac{2^{\frac{m}{m-1}}}{t_0}\int_{\Omega}u(t_0,x)G_{\Omega}(x , x_0)\dx
\le c_{1,\Omega_0} \frac{2^{\frac{m}{m-1}}}{t_0}\int_{B_r(x_0)\cup\Omega_r}\frac{u(t_0,x)}{|x-x_0|^{N-2s}}\dx\\
&\le c_{1,\Omega_0}\frac{2^{\frac{m}{m-1}}}{t_0}
    \left[\|u(t_0)\|_{\LL^\infty(\Omega)}\int_{B_r(x_0)}\frac{1}{|x-x_0|^{N-2s}}\dx
        +\int_{\Omega_r}\frac{u(t_0,x)}{|x-x_0|^{N-2s}}\dx\right]\\
\end{split}\end{equation*}
Next, we use the Young inequality $\alpha\beta\le \frac{1}{2}\alpha^{m}+2^{\frac{1}{m-1}}\beta^{\frac{m}{m-1}}$\,, valid for all $\alpha,\beta\ge 0$ and all $m >1$, together with
\begin{equation}\label{Upper.PME.Step.1.0}
 \int_{B_r(x_0)}\frac{1}{|x-x_0|^{N-2s}}\dx =\frac{\omega_N}{2s} r^{2s}\,,
\end{equation}
to obtain
\begin{equation*}\begin{split}
u^m (t_0,x_0) &\le \frac{1}{2}\|u(t_0)\|_{\LL^\infty(\Omega)}^m + 2^{\frac{1}{m-1}}\left[c_{1,\Omega_0}\frac{2^{\frac{m}{m-1}}}{t_0}
    \frac{\omega_N}{2s} r^{2s}\right]^{\frac{m}{m-1}}
        +c_{1,\Omega_0}\frac{2^{\frac{m}{m-1}}}{t_0\,r^{N-2s}}\int_{\Omega}u(t_0,x)\dx \\
&:=\frac{1}{2}\|u(t_0)\|_{\LL^\infty(\Omega)}^m +\left[A\, r^{\frac{2s m}{m-1}}+\frac{B}{r^{N-2s}}\right]:=\frac{1}{2}\|u(t_0)\|_{\LL^\infty(\Omega)}+F(r)\\
\end{split}
\end{equation*}
Therefore, we have proved that
\begin{equation}\label{Upper.PME.Step.1.1.a}
\|u(t_0)\|_{\LL^\infty(\Omega)} \le 2 F(r)=2\left[A\, r^{\frac{2s m}{m-1}}+\frac{B}{r^{N-2s}}\right]\,,
\end{equation}
where
\begin{equation}\label{Upper.PME.Step.1.2.a}
A:= 2^{\frac{1}{m-1}}\left(c_{1,\Omega}\frac{2^{\frac{m}{m-1}}}{t_0}
    \frac{\omega_N}{2s}\right)^{\frac{m}{m-1}}\qquad\mbox{and}\qquad
B:=c_{1,\Omega}\frac{2^{\frac{m}{m-1}}}{t_0}\|u(t_0)\|_{\LL^1(\Omega)}
\end{equation}
Then we choose $r=(B/A)^{(m-1)\vartheta_1}$
so that
\[\begin{split}
2F(r)&= 4A^{(N-2s)(m-1)\vartheta_1}B^{2sm\vartheta_1}
:=K_3 \frac{\|u(t_0)\|_{\LL^1(\Omega)}^{2sm\vartheta_1}}{t_0^{Nm\vartheta_1}}
\end{split}
\]
where we have defined $K_3$
\begin{equation}\label{Upper.PME.Step.1.3.a}
K_3:=2^{(N-2s)\vartheta_1+2}\left(\frac{\omega_N}{2s}
    \right)^{m(N-2s)\vartheta_1}
    \left(c_{1,\Omega}2^{\frac{m}{m-1}}  \right)^{mN\vartheta_1}
\end{equation}
This concludes the proof of the smoothing effect estimate \eqref{thm.Upper.PME.Smoothing.1}\,, once we recall that the  $\LL^1$-norm $\int_{\Omega} u(t_0,x)\dx$ is monotonically decreasing in time, cf. \cite{DPQRV1,DPQRV2}\,.\normalcolor

\noindent$\bullet~$\textit{$\LL^1_{\Phi_1}-\LL^\infty$ Smoothing effects via Type II estimates. } We will use the fundamental upper estimates \eqref{Upper.PME.Step.1.2} of Step 1 of the proof of Theorem \ref{thm.Upper.PME}\,, together with Type II estimates \eqref{typeII.Green.est}\,, as follows
\begin{equation*}\begin{split}
u^m (t_0,x_0) &\le \frac{2^{\frac{m}{m-1}}}{t_0}\int_{\Omega}u(t_0,x)G_{\Omega}(x , x_0)\dx\\
&\le \frac{2^{\frac{m}{m-1}}}{t_0}\left[\int_{B_r(x_0)}u(t_0,x)G_{\Omega}(x , x_0)\dx+\int_{\Omega_r}u(t_0,x)G_{\Omega}(x , x_0)\dx\right]\\
~_{(a)}&\le c_{1,\Omega_0}\frac{2^{\frac{m}{m-1}}}{t_0}
    \left[\|u(t_0)\|_{\LL^\infty(\Omega)}\int_{B_r(x_0)}\frac{1}{|x-x_0|^{N-2s}}\dx
        +\int_{\Omega_r}\frac{u(t_0,x)\Phi_1(x)}{|x-x_0|^{N-2s+1}}\dx\right]\\
~_{(b)}&\le \frac{1}{2}\|u(t_0)\|_{\LL^\infty(\Omega)}^m + 2^{\frac{1}{m-1}}\left[c_{1,\Omega_0}\frac{2^{\frac{m}{m-1}}}{t_0}
    \frac{\omega_N}{2s} r^{2s}\right]^{\frac{m}{m-1}}
        +c_{1,\Omega_0}\frac{2^{\frac{m}{m-1}}}{t_0\,r^{N-2s+1}}\int_{\Omega}u(t_0,x)\Phi_1(x)\dx \\
&:=\frac{1}{2}\|u(t_0)\|_{\LL^\infty(\Omega)}^m +\left[A\, r^{\frac{2s m}{m-1}}+\frac{B}{r^{N-2s+1}}\right]:=\frac{1}{2}\|u(t_0)\|_{\LL^\infty(\Omega)}+F(r)\\
\end{split}
\end{equation*}
where in $(a)$ we have used  the Green function estimates \eqref{typeII.Green.est}\,, namely $G_{\Omega}(x,x_0)\le c_{1,\Omega}|x-x_0|^{-(N-2s)}$ on the ball $B_r(x_0)$\,, while we have used  namely $G_{\Omega}(x,x_0)\le c_{1,\Omega}\Phi_1(x)|x-x_0|^{-(N-2s+1)}$ on $\Omega_r$\,. In $(b)$  we have used the Young inequality $\alpha\beta\le \frac{1}{2}\alpha^{m}+2^{\frac{1}{m-1}}\beta^{\frac{m}{m-1}}$\,, valid for all $\alpha,\beta\ge 0$ and all $m >1$; we have also used \eqref{Upper.PME.Step.1.0}\,. Therefore, we have proved that
\begin{equation}\label{Upper.PME.Step.1.1.b}
\|u(t_0)\|_{\LL^\infty(\Omega)} \le 2 F(r)=2\left[A\, r^{\frac{2s m}{m-1}}+\frac{B}{r^{N-2s+1}}\right]\,,
\end{equation}
where
\begin{equation}\label{Upper.PME.Step.1.2.b}
A:= 2^{\frac{1}{m-1}}\left(c_{1,\Omega}\frac{2^{\frac{m}{m-1}}}{t_0}
    \frac{\omega_N}{2s}\right)^{\frac{m}{m-1}}\qquad\mbox{and}\qquad
B:=c_{1,\Omega}\frac{2^{\frac{m}{m-1}}}{t_0}\|u(t_0)\|_{\LL^1_{\Phi_1}(\Omega)}
\end{equation}
Then we can choose $r=(B/A)^{(m-1)\vartheta_{1,1}}$, where we recall that $\vartheta_{1,1}=1/[2s+(N+1)(m-1)]$\,,
and that $1-2sm\vartheta_{1,1}=(m-1)(N-2s+1)\vartheta_{1,1}$\,, so that we obtain
\[\begin{split}
2F(r)&= 4A^{(N-2s+1)(m-1)\vartheta_{1,1}}B^{2sm\vartheta_{1,1}}
:=K_4^m \frac{\|u(t_0)\|_{\LL^1_{\Phi_1}(\Omega)}^{2sm\vartheta_{1,1}}}{t_0^{(N+1)m\vartheta_{1,1}}}
\end{split}
\]
where $K_4$ has the following expression
\begin{equation}\label{Upper.PME.Step.1.3.b}
K_4^m:=2^{(N-2s+1)\vartheta_{1,1}+2}\left(\frac{\omega_N}{2s}
    \right)^{m(N-2s+1)\vartheta_{1,1}}
    \left(c_{1,\Omega}2^{\frac{m}{m-1}}  \right)^{m(N+1)\vartheta_{1,1}}
\end{equation}
We conclude by observing that the weighted $\LL^1$-norm $\int_{\Omega} u(t_0,x)\Phi_1(x)\dx$ is monotone decreasing in time (see \eqref{L1weight.PME.Step.1.2} for a detailed explanation), so that \eqref{thm.Upper.PME.Smoothing.1} follows.\qed

\medskip

\noindent\textbf{Remark. }We shall observe that in order to get estimates for $u^m \in \LL^1\left((0,\infty):\LL^1_{\Phi_1}(\Omega)\right)$, a space required by our concept of solution, we need both the absolute bound and the smoothing effect of Theorems \eqref{thm.Upper.PME} and \eqref{thm.Upper.1.PME} respectively, indeed
\begin{equation}\label{remark.upper.final}\begin{split}
0\le & \int_0^{\infty}\int_{\Omega} u^m(t,x)\Phi_1(x)\dx\dt
\le \int_0^1\int_{\Omega} u^m(t,x)\Phi_1(x)\dx\dt +\int_1^{\infty}\int_{\Omega} u^m(t,x)\Phi_1(x)\dx\dt \\
&\le \int_0^1\|u(t)\|_{\LL^\infty(\Omega)}^{m-1}\int_{\Omega} u(t,x)\Phi_1(x)\dx\dt +\int_1^{\infty}\|u(t)\|_{\LL^\infty(\Omega)}^m\int_{\Omega}\Phi_1(x)\dx\dt \\
&\le \int_0^1\frac{K_4^{m-1}}{t^{(N+1)(m-1)\vartheta_{1,1}}}\|u_0\|_{\LL^1_{\Phi_1}(\Omega)}^{2s(m-1)\vartheta_{1,1}+1}\dt +\int_1^{\infty}\frac{K_1^m\|\Phi_1\|_{\LL^1(\Omega)}}{t^{\frac{m}{m-1}}}\dt <+\infty\\
\end{split}\end{equation}
the first integral is finite since $(N+1)(m-1)\vartheta_{1,1}<1$ and the second since $m/(m-1)>1$\,. We have also used the monotonicity in $t$ of $\|u(t)\|_{\LL^1_{\Phi_1}(\Omega)}$. \normalcolor

As a corollary of the above Theorem \ref{thm.Upper.1.PME}, we get the following reverse in time smoothing effects.
\begin{cor}[Backward Smoothing effects]\label{thm.Upper.Backward.PME}
There exists a universal constant $K_4>0$ such that for all $t,h>0$
\begin{equation}\label{thm.Upper.Backward.PME.2}
\| u(t)\|_{\LL^\infty(\Omega)} \le \frac{K_4}{t^{(d+1)\vartheta_{1,1}}}\left(1\vee \frac{h}{t}\right)^{\frac{2s\vartheta_{1,1}}{m-1}}\|u(t+h)\|_{\LL^1_{\Phi_1}(\Omega)}^{2s\vartheta_{1,1}}\mbox{\,.}
\end{equation}
\end{cor}
\noindent {\sl Proof.~}By the monotonicity estimates of Lemma \ref{lem.BC.est}\,, the function $u(x,t)t^{1/(m-1)}$ is non-decreasing in time for fixed $x$, therefore using the smoothing effect \eqref{thm.Upper.PME.Smoothing.2} we get for all $t_1\ge t$:
\begin{equation*}
\|u(t)\|_{\LL^\infty(\Omega)}
    \le \frac{K_4}{t^{(N+1)\vartheta_{1,1}}}\left(\int_\Omega u(t,x)\Phi_1(x)\dx\right)^{2s\vartheta_{1,1}}
    \le \frac{K_4}{t^{(N+1)\vartheta_{1,1}}}\left(\frac{t_1^{\frac{1}{m-1}}}{t^{\frac{1}{m-1}}}\int_\Omega u(t_1,x)\Phi_1(x)\dx\right)^{2s\vartheta_{1,1}}
\end{equation*}
where $K_4$ is given in Theorem \ref{thm.Upper.1.PME}\,. The above inequality implies \eqref{thm.Upper.Backward.PME.2} letting $t_1=t+h$\,.\qed

\section{Weighted $\LL^1$-estimates}
\label{Sect.Backward}

As an interesting application of the upper estimates of the previous section, we obtain the following weighted estimates, which will be very useful in the proof of the lower bounds in the next section.

\begin{prop}[Weighted $\LL^1$-estimates]\label{thm.L1weight.PME}
Under the current assumptions on $m$ and $u$, the integral $\int_{\Omega}u(t,x)\Phi_1(x)\dx$ is monotonically non-increasing in time and for all $0\le \tau_0\le \tau,t ~<+\infty$ we have
\begin{equation}\label{L1weight.PME.estimates.1}
\int_{\Omega}u(\tau,x)\Phi_1(x)\dx\le \int_{\Omega}u(t,x)\Phi_1(x)\dx + K_5\,\left|t-\tau\right|^{2s\vartheta_{1,1}}\left(\int_{\Omega} u(\tau_0)\Phi_1\dx\right)^{2s(m-1)\vartheta_{1,1}+1}
\end{equation}
where $K_5:=\,\lambda_1\,K_4 /(2s\vartheta_{1,1})$ and $K_4>0$ is given in Theorem $\ref{thm.Upper.1.PME}$\,.
\end{prop}

\noindent\textsl{Proof of Proposition \ref{thm.L1weight.PME}. }Let us (formally) calculate
\begin{equation}\label{L1weight.PME.Step.1.1}
\begin{split}
\frac{\rd}{\dt}\int_{\Omega}u(t,x)\Phi_1(x)\dx
&=-\int_{\Omega}\A (u^m)\Phi_1\dx
=-\int_{\Omega}u^m\A \Phi_1\dx
=-\lambda_1\int_{\Omega}u^m\Phi_1\dx\le 0
\end{split}
\end{equation}
where we have used Definition \eqref{Def.Very.Weak.Sol.Dual} of very weak solution, together with the fact that $\A\Phi_1=\lambda_1\Phi_1\ge 0$\,. Integrating the above inequality we obtain the monotonicity property for all $0\le \tau_0\le t$:
\begin{equation}\label{L1weight.PME.Step.1.21}
\int_{\Omega}u(t,x)\Phi_1(x)\dx=\int_{\Omega}u(\tau_0,x)\Phi_1(x)\dx-\lambda_1\int_{\tau_0}^{t}\int_{\Omega}u^m(\tau,x)\Phi_1\dx \rd \tau
\le \int_{\Omega}u(\tau_0,x)\Phi_1(x)\dx\,.
\end{equation}
We will obtain better bounds in what follows. We just remark that the derivation of the above equality \eqref{L1weight.PME.Step.1.1} made here is formal: the correct form is the equality \eqref{L1weight.PME.Step.1.21}, which has already been proved  in Step 1 of the proof of Proposition \ref{thm.NLE.PME}, more precisely, equality \eqref{L1weight.PME.Step.1.21} is exactly equality \eqref{step.1.thm.repr} with the (admissible) choice of the test function $\psi=\Phi_1$\,. From inequality \eqref{L1weight.PME.Step.1.21} it immediately follows the monotonicity property
\begin{equation}\label{L1weight.PME.Step.1.2}
\int_{\Omega}u(t,x)\Phi_1(x)\dx\le\int_{\Omega}u(\tau,x)\Phi_1(x)\dx\qquad\mbox{for all }0\le \tau\le t\,.
\end{equation}
On the other hand, using the weighted smoothing effect \eqref{thm.Upper.PME.Smoothing.2} and \eqref{L1weight.PME.Step.1.2}\,, we get for all $t\ge   \tau_0$\,:
\begin{equation}\label{L1weight.PME.Step.1.3}
\begin{split}
\lambda_1\int_{\tau_0}^t\int_{\Omega}u^m(\tau)\Phi_1\dx\rd\tau
&\le \lambda_1\,\int_{\tau_0}^t \|u(t)\|_{\LL^\infty(\Omega)}^{m-1}\, \int_{\Omega}u(t)\Phi_1\dx\rd\tau\\
&\le  \int_{\tau_0}^t\frac{\lambda_1\,K_4}{(t-\tau_0)^{(N+1)(m-1)\vartheta_{1,1}}}\rd\tau\left(\int_{\Omega} u(\tau_0)\Phi_1\dx\right)^{2s(m-1)\vartheta_{1,1}+1}\\
\end{split}
\end{equation}
Recalling that $1-2sm\vartheta_{1,1}=(m-1)(N-2s+1)\vartheta_{1,1}$, we combine \eqref{L1weight.PME.Step.1.21} with \eqref{L1weight.PME.Step.1.3} to get for all $t,\tau\ge \tau_0$:
\[
\begin{split}
\int_{\Omega}u(\tau)\Phi_1\dx\le \int_{\Omega}u(t)\Phi_1\dx
&+ \frac{\lambda_1\,K_4}{2s\vartheta_{1,1}}\,
    \left(\int_{\Omega} u(\tau_0)\Phi_1\dx\right)^{2s(m-1)\vartheta_{1,1}+1}\left|(t-\tau_0)^{2s\vartheta_{1,1}}-(\tau-\tau_0)^{2s\vartheta_{1,1}}\right|\,,
\end{split}
\]
which implies \eqref{L1weight.PME.estimates.1} using the numerical inequality $|(t-\tau_0)^{2s\vartheta_{1,1}}-(\tau-\tau_0)^{2s\vartheta_{1,1}}|\le |t-\tau|^{2s\vartheta_{1,1}}$\,, valid since $2s\vartheta_{1,1}<1$.\qed

The above Proposition has interesting consequences, the first one is the following Corollary which is crucial in the proof of the lower bounds of the next section.

\begin{cor}[Backward in time $\LL^1_{\Phi_1}$ lower bounds]\label{cor.abs.L1.phi}  For all
\begin{equation}\label{L1weight.PME.estimates.2}
0\le \tau_0\le  t\le \tau_0+ \frac{1}{(2K_5)^{1/(2s\vartheta_{1,1})}\left(\int_{\Omega} u(\tau_0)\Phi_1\dx\right)^{m-1}}
\end{equation}
we have
\begin{equation}\label{L1weight.PME.estimates.3}
\frac{1}{2}\int_{\Omega}u(\tau_0,x)\Phi_1(x)\dx\le \int_{\Omega}u(t,x)\Phi_1(x)\dx\,.
\end{equation}
where $K_5>0$ is as in Proposition $\ref{thm.L1weight.PME}$\,.
\end{cor}
\noindent\textsl{Proof of Corollary \ref{cor.abs.L1.phi}. }The proof of \eqref{L1weight.PME.estimates.3} follows from \eqref{L1weight.PME.estimates.1} by letting $\tau=\tau_0$ and by choosing $|t-\tau_0|$ ``small'', namely as in \eqref{L1weight.PME.estimates.2}\,.\qed

\medskip

\normalcolor


\section{Quantitative positivity estimates}
\label{sect.positivity}
We use upper pointwise estimates of Theorem \ref{thm.NLE.PME} to prove lower bounds in the form of lower Harnack inequalities. We will also use the upper estimates proved in the previous section together with the $\LL^1_{\Phi_1}(\Omega)$-estimates\,. The following lower estimates are consequence of Type II estimates for the Green function \eqref{typeII.Green.est}.
\begin{thm}[Lower absolute and boundary estimates]\label{thm.Lower.PME}
Let let $m>1$ and let $u\ge 0$ be a solution in the class $\mathcal{S}$ of very weak solutions to the Dirichlet problem \eqref{FPME.prob}, corresponding to the initial datum $0\le u_0\in \LL^1_{\Phi_1}(\Omega)$\,. Then, there exist  constants $L_0(\Omega),L_1(\Omega)>0$\,, so that, setting
\begin{equation}\label{thm.Lower.PME.Boundary.0}
t_*= \frac{L_0(\Omega)}{\left(\int_{\Omega}u_0\Phi_1\dx\right)^{m-1}}\,,
\end{equation}
we have that for all $t\ge t_*$ and all $x_0\in \Omega$, the following inequality holds:
\begin{equation}\label{thm.Lower.PME.Boundary.1}
u(t,x_0)\ge L_1(\Omega)\,\frac{\Phi_1(x_0)^{\frac{1}{m}}}{t^{\frac{1}{m-1}}}\,.
\end{equation}
\end{thm}
The constants $L_0(\Omega),L_1(\Omega)>0$\,, depend on $N, m, s$ and on $\Omega$\,, but not on $u$ (or any norm of $u$); they have an explicit form given at the end of the proof.

\noindent {\sl Proof of Theorem \ref{thm.Lower.PME}.~} \textsc{Step 1. }\textit{Fundamental lower estimates. }We first recall the upper pointwise estimates of Theorem \ref{thm.NLE.PME}, namely for all $0\le t_0\le t_1 \le t $ and $x_0\in \Omega$\,, we have that
\begin{equation}\label{Lower.PME.Step.1.1}
\int_{\Omega}u(t_0,x)G_{\Omega}(x , x_0)\dx - \int_{\Omega}u({t_1},x)G_{\Omega}(x , x_0)\dx \le \frac{m-1}{t_0^{\frac{1}{m-1}}} t^{\frac{m}{m-1}}\,u^m(t,x_0)\,,
\end{equation}
Then we estimate:
\begin{equation}\label{Lower.PME.Step.1.2}\begin{split}
\int_{\Omega}u({t_1},x)G_{\Omega}(x , x_0)\dx
&\le_{(a)}\overline{K}_2 \frac{\Phi_1(x_0)}{t_1^{\frac{1}{m-1}}}
\le_{(b)}\frac{1}{2}\int_{\Omega}u(t_0,x)G_{\Omega}(x , x_0)\dx
\end{split}
\end{equation}
where in $(a)$ we have used the absolute upper bounds \eqref{thm.Upper.PME.Boundary.2} of the remark after Theorem \eqref{thm.Upper.PME.Boundary}, namely for  all $t_1> 0$ and $x_0\in \Omega$ we have $\int_{\Omega}u(t_1,x)G_{\Omega}(x , x_0)\dx\le \overline{K}_2 \Phi_1(x_0)t_1^{-\frac{1}{m-1}}$\,, where $\overline{K}_2=c_{5,\Omega}^m 2^{\frac{m}{(m-1)^2}}$ with $c_{5,\Omega}>0$ is the universal constant given in Lemma \ref{Lem.Green.2}\,. In $(b)$ we just have chosen $t_1$ ``relatively big'', namely
\begin{equation}\label{Lower.PME.Step.1.3}
t_1\ge \left[\frac{2\,\overline{K}_2\,\Phi_1(x_0)}{\int_{\Omega}u(t_0,x)G_{\Omega}(x , x_0)\dx }\right]^{m-1}\,,
\end{equation}
where we have used that the weighted $\LL^1$-norm $\int_{\Omega}u(t_0,x)G_{\Omega}(x , x_0)\dx$ is non increasing in time.
Joining the above inequalities, we obtain that for all $0\le t_0\le t_1 \le t $
\begin{equation}\label{Lower.PME.Step.1.4}
\frac{t_0^{\frac{1}{m-1}}}{2(m-1)}\int_{\Omega}u(t_0,x)G_{\Omega}(x , x_0)\dx \le  t^{\frac{m}{m-1}}\,u^m(t,x_0)\,.
\end{equation}
Finally, we notice that we always have $t_0\le t_1$\,, thanks to  the absolute upper bounds \eqref{thm.Upper.PME.Boundary.2}\,, namely for all namely for  all $t_0> 0$ and $x_0\in \Omega$ we have $\int_{\Omega}u(t_0,x)G_{\Omega}(x , x_0)\dx\le \overline{K}_2 \Phi_1(x_0)t_0^{-\frac{1}{m-1}}$ so that
\begin{equation}\label{Lower.PME.Step.1.5}
t_0^{\frac{1}{m-1}}\le \frac{\overline{K}_2 \Phi_1(x_0)}{\int_{\Omega}u(t_0,x)G_{\Omega}(x , x_0)}
\le\frac{2\,\overline{K}_2 \Phi_1(x_0)}{\int_{\Omega}u(t_0,x)G_{\Omega}(x , x_0)}\le t_1^{\frac{1}{m-1}}
\end{equation}
for the moment we still keep $t_0$ and $t_1$ free, the only restriction is the lower bound \eqref{Lower.PME.Step.1.3} on $t_1$.

\noindent$\bullet~$\textsc{Step 2. }\textit{Absolute lower and boundary estimates. }It is sufficient now to apply the Type II lower Green function estimates \eqref{typeII.Green.est}, namely $c_{0,\Omega}\Phi_1(x)\Phi_1(x_0)\le G_{\Omega}(x,x_0)$\,, to the fundamental lower bound \eqref{Lower.PME.Step.1.4} to get that for all $0\le t_0\le t_1 \le t $
\begin{equation}\label{Lower.PME.Step.2.1}
\frac{c_{0,\Omega_0}\Phi_1(x_0)}{2(m-1)}t_0^{\frac{1}{m-1}}\int_{\Omega}u(t_0)\Phi_1\dx
\le \frac{t_0^{\frac{1}{m-1}}}{2(m-1)}\int_{\Omega}u(t_0)G_{\Omega}(\cdot , x_0)\dx \le  t^{\frac{m}{m-1}}\,u^m(t,x_0)\,,
\end{equation}
where $t_1$ must satisfy inequality \eqref{Lower.PME.Step.1.3}\,. Next we use the backward in time lower estimates of Corollary \ref{cor.abs.L1.phi}, with the particular choice $t_0= (2K_5)^{-1/(2s\vartheta_{1,1})}\left(\int_{\Omega} u(\tau_0)\Phi_1\dx\right)^{-(m-1)}
$ so that, inequality \eqref{L1weight.PME.estimates.3} implies that for such $t$ we have an absolute lower bound for the quantity:
\begin{equation}\label{L1weight.PME.estimates.4}\begin{split}
t_0^{\frac{1}{m-1}}
&\int_{\Omega}u(t_0)\Phi_1\dx
 \ge \frac{t_0^{\frac{1}{m-1}}}{2}\int_{\Omega}u_0\Phi_1\dx
=\frac{1}{2}(2K_5)^{\frac{-1}{(2s\vartheta_{1,1})(m-1)}}\,.
\end{split}
\end{equation}
Therefore we have that:
\[
t_0=  \frac{1}{(2K_5)^{1/(2s\vartheta_{1,\gamma})}\left(\int_{\Omega} u_0\Phi_1\dx\right)^{m-1}}
\qquad\mbox{implies}\qquad
t_0^{\frac{1}{m-1}}\int_{\Omega}u(t_0)\Phi_1\dx\ge\frac{1}{2}(2K_5)^{\frac{-1}{(2s\vartheta_{1,\gamma})(m-1)}}\,,
\]
and gives the following absolute lower bound for all $0\le t_0\le t_1 \le t $
\begin{equation}\label{Lower.PME.Step.2.2}
L_1^m\Phi_1(x_0):=\frac{c_{0,\Omega_0}\Phi_1(x_0)}{4(m-1)(2K_5)^{\frac{1}{2s\vartheta_{1,\gamma}(m-1)}}}\le \frac{c_{0,\Omega_0}\Phi_1(x_0)}{2(m-1)}t_0^{\frac{1}{m-1}}\int_{\Omega}u(t_0)\Phi_1\dx\le  t^{\frac{m}{m-1}}\,u^m(t,x_0)\,,
\end{equation}
where $K_5$ is given in Theorem \ref{thm.L1weight.PME}\,, and $t_0\le t_1$\,, namely inequality \eqref{Lower.PME.Step.1.5} continues to hold in the form:
\begin{equation}\label{Lower.PME.Step.2.3}
t_0:=\frac{1}{(2K_5)^{1/(2s\vartheta_{1,\gamma})}\left(\int_{\Omega} u_0\Phi_1\dx\right)^{m-1}}
\le \left[\frac{2\,\overline{K}_2 \Phi_1(x_0)}{\int_{\Omega}u(t_0,x)G_{\Omega}(x , x_0)\dx }\right]^{m-1}:= t_1\,.
\end{equation}
notice that we have fixed the value of $t_1$. The only thing which is not yet explicit is $t_1$\,, that depends on the norm of $u(t_0)$; in the next step we are going to show that how we can estimate from above such quantity in terms of the initial datum $u_0$\,.

\noindent$\bullet~$\textsc{Step 3. }\textit{The critical time $t_*$. }We want to find an upper bound for $t_1$, which will be the critical time $t_*$, which shall only depend on $u_0$ but not on $x_0$, namely
\begin{equation}\label{Lower.PME.Step.3.1}\begin{split}
t_1^{\frac{1}{m-1}} &=\frac{2\,\overline{K}_2 \Phi_1(x_0)}{\int_{\Omega}u(t_0,x)G_{\Omega}(x , x_0)\dx }
  \le_{(a)}\frac{4\,\overline{K}_2 \Phi_1(x_0)}{c_{0,\Omega}\Phi_1(x_0)\int_{\Omega}u_0\Phi_1\dx} \\
 &=\frac{4\,\overline{K}_2  }{c_{0,\Omega} \int_{\Omega}u_0\Phi_1\dx}
  := t_*^{\frac{1}{m-1}}\,,
\end{split}
\end{equation}
where in $(a)$ we have used again the backward in time lower estimates of Corollary \ref{cor.abs.L1.phi}, since for all
\[
t_0\le \frac{1}{(2K_5)^{1/(2s\vartheta_{1,\gamma})}\left(\int_{\Omega} u_0\Phi_1\dx\right)^{m-1}}
\qquad\mbox{we have}\qquad
\int_{\Omega}u(t_0)\Phi_1\dx\ge\frac{1}{2}\int_{\Omega}u_0\Phi_1\dx\,,
\]
which, combined with the Type II lower Green function estimates \eqref{typeII.Green.est}, namely $c_{0,\Omega}\Phi_1(x)\Phi_1(x_0)\le G_{\Omega}(x,x_0)$\,, gives
\[
\int_{\Omega}u(t_0,x)G_{\Omega}(x , x_0)\dx
\ge c_{0,\Omega}\Phi_1(x_0)\int_{\Omega}u(t_0)\Phi_1\dx
\ge\frac{c_{0,\Omega}}{2}\Phi_1(x_0)\int_{\Omega}u_0\Phi_1\dx\,.
\]
Therefore we have proved that for all $t\ge t_*\ge t_1$ inequality \eqref{thm.Lower.PME.Boundary.1} holds. We conclude the proof by providing an explicit expression of the constants.

\noindent$\bullet~$\textsc{Step 4. }\textit{Expression for the constants. }We have set
\begin{equation}\label{Lower.PME.Step.4.1}
t_*= \frac{L_0}{\left(\int_{\Omega}u_0\Phi_1\dx\right)^{m-1}}\,,\quad\mbox{with}\quad L_0=\left[\frac{4\,\overline{K}_2}{c_{0,\Omega}}\right]^{m-1}\,,\quad\mbox{and}\quad
L_1^m:=\frac{c_{0,\Omega_0} }{4(m-1)(2K_5)^{\frac{1}{2s\vartheta_{1,\gamma}(m-1)}}}\,,
\end{equation}
where $\overline{K}_2$ is given in the remark after Theorem \ref{thm.Upper.PME}\,, $K_5$ is given in Proposition \ref{thm.L1weight.PME}, $c_{0,\Omega}$ is the constant in Type II estimates \eqref{typeII.Green.est}\,.\qed

\medskip

\noindent\textbf{Remarks. } (i) Recall that $\Phi_1$ is the first eigenfunction of $\A$ and satisfies estimates \eqref{Phi1.est}:
\[
\Phi_1(x)\asymp  \dist(x, \partial\Omega) \wedge 1 \qquad\mbox{for all $x\in \Omega$.}
\]
Therefore, the lower boundary behaviour of $u(t,\cdot)$ is:
\begin{equation}\label{Lem.Green.est.Upper.II.B.3}
u(t,x)\ge \, \frac{L_1}{t_0^{\frac{1}{m-1}}} \big(\dist(x_0,\partial\Omega)^{\frac{1}{m}}\wedge 1\big)\,,\qquad\qquad\mbox{for all $t_0\ge t_*\ge 0$ and $x_0\in \Omega$\,.}
\end{equation}
This boundary behaviour is sharp because we have already obtained upper bounds with matching powers of $\Phi_1$, cf. Theorem \ref{thm.Upper.PME.Boundary}, and also because the solutions obtained by separation of variables have the same boundary behaviour, as already explained in the remark after Theorem \ref{thm.Upper.PME.Boundary}\,.\\
\noindent(ii) The above positivity estimate \eqref{thm.Lower.PME.Boundary.1} can be stated as a lower bound in terms of the separation of variables solution $V(x)(t+h)^{-1/(m-1)}$ for some $h>0$ that depends on the particular solution\,: we have that  $u(t,x)\ge V(x)(t+h)^{-1/(m-1)}$ for all $t\ge t_*$\,. This has important consequences for the study of the asymptotic behaviour, cf. \cite{BSV2013}\,.\\
\noindent(iii) $t_*$ is an estimate the time that it takes to fill the hole: if $u_0$ is concentrated close to the border (leaves an hole in the middle of $\Omega$), then $\int_{\Omega}u_0\Phi_1\dx$ is small, therefore $t_*$ becomes very large, therefore it takes a lot of time to fill the hole.\\
\noindent(iv)  These estimates can also be rewritten in an equivalent way as Aronson-Caffarelli type estimates, in the spirit of the estimates first proved in \cite{ArCaff} for $s=1$ and $\Omega=\RR^N$, and generalized by us in \cite{BV-ADV, BV2012}. We can rephrase the lower estimates \eqref{thm.Lower.PME.Boundary.1} as follows:
\[
\mbox{either}\quad t\le t_*= \frac{L_0}{\left(\int_{\Omega}u_0\Phi_1\dx\right)^{m-1}}\,,
\quad\mbox{or}\quad
u(t,x_0)\ge L_1\,\frac{\Phi_1(x_0)^{\frac{1}{m}}}{t^{\frac{1}{m-1}}}\quad\mbox{for all $t\ge t_*$ and all $x_0\in \Omega$}\,,
\]
which gives, for all $t\ge 0$ and all $x_0\in \Omega$:
\[
1\le \left(\frac{t_*}{t}\right)^{\frac{1}{m-1}}+\frac{t^{\frac{1}{m-1}}u(t,x_0)}{L_1\Phi_1(x_0)^{\frac{1}{m}}}\,,\qquad\mbox{that is}\qquad
 u(t,x_0)\ge\frac{L_1\Phi_1(x_0)^{\frac{1}{m}}}{t^{\frac{1}{m-1}}}\left[1-\left(\frac{t_*}{t}\right)^{\frac{1}{m-1}}\right]\,.
\]

As an open problem, it would be interesting to find precise lower bounds for small times, namely $0<t<t_*$.

\section{Harnack inequalities}
\label{sect.Haranck}
Let let $m>1$ and let $u\ge 0$ be a solution in the class $\mathcal{S}$ of very weak solutions to the Dirichlet problem \eqref{FPME.prob}, corresponding to the initial datum $0\le u_0\in \LL^1_{\Phi_1}(\Omega)$\,. We prove the following weak Harnack principle.

\begin{thm}\label{thm.GHP.PME}
There exist universal constants $H_0, H_1, L_0>0$ such that setting
\begin{equation}\label{thm.GHP.PME.0}
t_*= \frac{L_0}{\left(\int_{\Omega}u_0\Phi_1\dx\right)^{m-1}}\,,
\end{equation}
we  have that for all $t\ge t_*$ and all $x\in \Omega$, the following inequality holds:
\begin{equation}\label{thm.GHP.PME.1}
H_0\,\frac{\Phi_1(x)^{\frac{1}{m}}}{t^{\frac{1}{m-1}}} \le \,u(t,x)\le H_1\, \frac{\Phi_1(x)^{\frac{1}{m}}}{t^{\frac{1}{m-1}}}
\end{equation}
Recall that $\Phi_1$ is the first eigenfunction of $\A$ and satisfies estimates \eqref{Phi1.est}.
\end{thm}
\noindent {\sl Proof.~}We combine the upper bounds \eqref{thm.Upper.PME.Boundary} of Theorem \ref{thm.Upper.PME} with the lower bounds \eqref{thm.Lower.PME.Boundary.1} of Theorem \ref{thm.Lower.PME}; the expression of $t_*$ is explicitly given in Theorem \ref{thm.Lower.PME}\,, and the constant $H_0=L_1$\,, where $L_1$ is as in Theorem \ref{thm.Lower.PME}\,; as for the constant $H_1=K_2$\,, where $K_2$ is given in Theorem \ref{thm.Upper.PME}\,. Recall that constants $L_0, H_1,H_0>0$ depend only on $N, m, s,\gamma$ and $\Omega$\,, but not on $u$\,, and have an explicit form given in the proof. \qed

\medskip

These estimates are very useful for the asymptotic behaviour that we will study in a forthcoming paper, \cite{BSV2013}. \color{blue} In this nonlinear setting they do not directly imply $C^\alpha$ regularity, as in \cite{AC}. \color{darkblue}
On the other hand, we can show that the solution $u$ to the parabolic problem somehow inherits the (local) Harnack inequality that holds for the first eigenfunction $\Phi_1$ of the operator $\A$, namely for all $B_R(x_0)\in \Omega$\,:
\begin{equation}\label{Harnack.Phi1}
\sup_{x\in B_R(x_0)}\Phi_1(x)\le \mathcal{H}\inf_{x\in B_R(x_0)}\Phi_1(x)
\end{equation}
where the constant $\mathcal{H}>0$ is universal and can be made explicit, see for example \cite{BGV-Milan}. \normalcolor

\begin{thm}[Local Harnack Inequalities of Elliptic Type]\label{thm.Harnack.Local}
There exist universal constants $H_0$, $H_1$, $L_0>0$ such that setting
\begin{equation}\label{thm.Harnack.PME.0}
t_*= \frac{L_0}{\left(\int_{\Omega}u_0\Phi_1\dx\right)^{m-1}}\,,
\end{equation}
we  have that for all $t\ge t_*$ and all $B_R(x_0)\in \Omega$, the following inequality holds:
\begin{equation}\label{thm.Harnack.PME.1}
\sup_{x\in B_R(x_0)}u(t,x)\le \frac{H_1\,\mathcal{H}^\frac{1}{m}}{H_0}\, \inf_{x\in B_R(x_0)}u(t,x)
\end{equation}
\end{thm}
\noindent {\sl Proof.~}\color{darkblue}We combine the above weak Harnack inequality \eqref{thm.GHP.PME.1} with the Harnack inequality \eqref{Harnack.Phi1} that holds for $\Phi_1$ and we obtain \eqref{thm.Harnack.PME.1}:
\begin{equation}\label{Harnack.PME.step.1.1}
\sup_{x\in B_R(x_0)}u(t,x)
\le \frac{H_1}{t^{\frac{1}{m-1}}}\sup_{x\in B_R(x_0)}\Phi_1^\frac{1}{m}(x)
\le \frac{H_1\,\mathcal{H}^\frac{1}{m}}{t^{\frac{1}{m-1}}}\inf_{x\in B_R(x_0)}\Phi_1^\frac{1}{m}(x)
\le \frac{H_1\,\mathcal{H}^\frac{1}{m}}{H_0}\inf_{x\in B_R(x_0)}u(t,x)\,,
\end{equation}
where we recall that $\mathcal{H}>0$ is the constant in \eqref{Harnack.Phi1}\,.\qed\normalcolor

\medskip

\noindent\textbf{Remarks. } (i) Also in the case $s=1$\,, these Harnack inequalities are stronger than the known Harnack inequalities  for the porous medium equation, , cf. \cite{D88, DGVbook}, which are of forward type and are stated in terms of the so-called intrinsic geometry. On the other hand, elliptic Harnack inequalities are false for $m=1$, as shown by simple examples using the fundamental solution.
Indeed, in the linear case only forward in time Harnack inequalities are possible, therefore we are observing a phenomenon typical of non-linear diffusion. Usually, Harnack inequalities of elliptic type occur in the fast diffusion range, namely when $m<1$\,, cf. our works \cite{BGV-Domains, BV, BV-ADV, BV2012}; the above estimates are surprisingly valid for $m>1$ and to our knowledge they have never been observed before.\normalcolor\\
(ii) Even if the operator is nonlocal, the reader may like to find a ``purely local'' Harnack estimate, i.\,e., one that holds depending on local data. Such inequality holds after a longer time that still depends on the global norm of $u_0$. Indeed we could have defined a ``bigger'' $\overline{t}_*$ in terms of a local norms as follows
\begin{equation}\label{thm.Harnack.PME.4}
t_*=\frac{L_0}{\left(\int_{\Omega}u_0\Phi_1\dx\right)^{m-1}}\le \frac{L_2}{\dist\big(B_R(x_0)\,,\partial\Omega\big)\left(\int_{B_R(x_0)}u_0\dx\right)^{m-1}}:=\overline{t}_*\,.
\end{equation}
and estimate \eqref{thm.Harnack.PME.1} holds after the local time  $\overline{t}_*$.


\section{Bounds and boundary behaviour of solution to the elliptic problem}
\label{sect.Elliptic.Prob}

We consider the homogeneous Dirichlet problem
\begin{equation}\label{Elliptic.prob}
\left\{\begin{array}{lll}
\A(V^m)= \lambda V\,, &  ~ {\rm in}~  \Omega\,,\\
V=0\,, & ~\mbox{on $\partial \Omega$\,,}
\end{array}
\right.
\end{equation}
where $\Omega\subset \RR^d$ is a bounded domain with smooth boundary, $\lambda>0$, $m>1$, $0<s\le 1$ and $d\ge 1$\,; the linear operator $\A$ will be either the spectral fractional Laplacian, $(-\Delta_{\Omega})^{s}$, or the regional  fractional Laplacian, $(-\Delta_{|\Omega})^{s}$, that we have previously introduced and discussed in Section \eqref{sect.prelim}\,.

\begin{defn}\label{Def.Very.Weak.Sol.Dual.elliptic} A function $V$ is a {\sl very weak} solution to Problem \eqref{Elliptic.prob} if:
\begin{itemize}
\item $V\in \LL^1_{\Phi_1}(\Omega)$\,, $|V|^{m-1}V \in \LL^1_{\Phi_1}(\Omega)$;
\item  The identity
\begin{equation}
\lambda\int_{\Omega}\AI (V) \,\psi\,\dx
=\int_{\Omega} |V|^{m-1}V\,\psi\,\dx.
\end{equation}
holds for every test function $\psi$ such that $\psi/\Phi_1\in \LL^\infty(\Omega)$\,.\normalcolor
\end{itemize}
\end{defn}
\noindent\textbf{Remark. }The above definition of very weak solution, has been given in \cite{KR07}, for $s=1$.\normalcolor
\begin{thm}[Bounds and boundary behaviour for the elliptic problem]\label{Thm.Elliptic.Harnack}
Let $V\ge~ 0$ be a very weak solution to the Dirichlet Problem \eqref{Elliptic.prob}, then there exist universal positive constants $h_0$ and $h_1$ such that the following estimates hold true for all $x_0\in \Omega$:
\begin{equation}\label{Thm.Elliptic.Harnack.ineq}
h_0\|V\|_{\LL^1_{\Phi_1}}\Phi_1(x_0)\le V^m(x_0)\le h_1\Phi_1(x_0)\,,
\end{equation}
where $h_1=c_{5,\Omega}\lambda^{1/(m-1)}$ and $h_0=c_{0,\Omega}\lambda$\,, with $c_{5,\Omega}$ given in Lemma $\ref{Lem.Green.2}$ and $c_{0,\Omega}$ is the constant in the Type II lower estimates \eqref{typeII.Green.est}.
\end{thm}
\noindent {\sl Proof.~}We split two steps.

\noindent$\bullet~$\textsc{Step 1. }\textit{A pointwise equality. }We want to prove that for almost all $x_0\in\Omega$ we have
\begin{equation}\label{Step.1.1Thm.Elliptic.Harnack.ineq}
V^m(x_0)=\int_{\Omega}V(x) \,G_\Omega(x_0,x)\,\dx\,.
\end{equation}
To prove this formula, we first use Definition \ref{Def.Very.Weak.Sol.Dual.elliptic} of very weak solution to get:
\begin{equation}\label{Step.1.2Thm.Elliptic.Harnack.ineq}
\int_{\Omega} V^m\,\psi\,\dx=\lambda\int_{\Omega}\AI (V) \,\psi\,\dx=\lambda\int_{\Omega}V \,\AI\psi\,\dx
\end{equation}
for any $\psi$ such that $\psi/\Phi_1\in \LL^\infty(\Omega)$\,.  The proof of formula \eqref{Step.1.1Thm.Elliptic.Harnack.ineq} now follows by approximating the Green function $G_\Omega(x_0,\cdot)$ by means of a sequence of admissible test functions $\psi_{n}^{(x_0)}$, as it has been done in Step 4 of the proof of Theorem \ref{thm.NLE.PME}\,.

\noindent$\bullet~$\textsc{Step 2. }\textit{The application of Integral Green Estimates I and II. }Inequality \eqref{Step.1.1Thm.Elliptic.Harnack.ineq} of Step 1, guarantees hypothesis \eqref{Lem.Green.2.hyp} with $\lambda=\kappa_0$.  Therefore we use the integral estimates II of Lemma \ref{Lem.Green.2} to obtain that
\[
V^m(x_0)\le c_{5,\Omega}^m \lambda^{\frac{m}{m-1}}\Phi_1(x_0)\,,
\]
which proves the upper bound of inequality \eqref{Thm.Elliptic.Harnack.ineq}\,. The lower bound of inequality \eqref{Thm.Elliptic.Harnack.ineq} directly follow by inequality \eqref{Lem.Green.est.Lower.II} of Lemma \ref{Lem.Green} or simply from Type II lower estimates \eqref{typeII.Green.est}\,, namely $G_\Omega(x_0,x)\ge c_{0,\Omega}\Phi_1(x)\Phi_1(x_0)$\,, to get
\[
V^m(x_0)= \lambda \int_{\Omega}V(x) \,G_\Omega(x_0,x)\,\dx \ge c_{0,\Omega} \lambda \|V\|_{\LL^1_{\Phi_1}(\Omega)} \Phi_1(x_0)\,.
\]
The proof of inequality \eqref{Thm.Elliptic.Harnack.ineq} is concluded.\qed
\normalcolor


\section{Existence and uniqueness of non-negative very weak solutions}\label{sect.Exist.Uniq}

We begin with a result on almost weighted $\LL^1$-contractivity for ordered solutions

\begin{prop}\label{thm.L1weight.contr}
 let $u\ge v$ be two ordered very weak solutions to the Dirichlet problem \eqref{FPME.prob} corresponding to the initial data $0\le u_0,v_0\in \LL^1_{\Phi_1}(\Omega)$\,. Then the integral $\int_{\Omega}\big[u(t,x)-v(t,x)\big]\Phi_1(x)\dx$ is monotonically non-increasing in time and for all $0\le \tau_0\le \tau,t ~<+\infty$ we have
\begin{equation}\label{L1weight.contr.estimates.1}\begin{split}
\int_{\Omega}\big[u(\tau,x)-v(\tau,x)\big]\Phi_1(x)\dx
&\le \int_{\Omega}\big[u(t,x)-v(t,x)\big]\Phi_1(x)\dx \\
&+ K_7[u(\tau_0),v(\tau_0)]\,\left|t-\tau\right|^{2s\vartheta_{1,1}}\,\int_{\Omega} \big[u(\tau_0,x)-v(\tau_0,x)\big]\Phi_1\dx
\end{split}
\end{equation}
where
\begin{equation}\label{L1weight.contr.estimates.1.a}
 K_7:= \frac{m K_4^{m-1} \lambda_1}{2s\vartheta_{1,1}}
     \max\left\{\|u(\tau_0)\|_{\LL^1_{\Phi_1}(\Omega)}^{2s\vartheta_{1,1}}
            \,,\,\|v(\tau_0)\|_{\LL^1_{\Phi_1}(\Omega)}^{2s\vartheta_{1,1}}\right\}^{m-1} \,
\end{equation}
and $K_4>0$ is given in Theorem $\ref{thm.Upper.1.PME}$\,.
\end{prop}
\noindent {\sl Proof.~}We begin by applying to $u$ and $v$ separately the definition of very weak solution \ref{Def.Very.Weak.Sol.Dual} in the form given in formula \eqref{step.1.thm.repr} of Step 1 of the proof of Proposition \ref{thm.NLE.PME}, with the admissible test function $\psi=\Phi_1$ (recall that $\AI \Phi_1=\lambda_1^{-1}\Phi_1\ge 0$); therefore we get for any $t,t_0\ge 0$
\begin{equation}\label{L1weight.contr.Step.1.1}
\begin{split}
\int_{\Omega}\big[u(t_0,x)-v(t_0,x)\big]\Phi_1(x)\dx- \int_{\Omega}\big[u(t,x)-v(t,x)\big]\Phi_1(x)\dx
&=\lambda_1\int_{t_0}^{t_1}\int_{\Omega} (u^m-v^m)\Phi_1\dx\rd\tau\\
\end{split}
\end{equation}
From this equality and the fact that solutions are ordered, namely $u-v\ge 0$, it immediately follows the monotonicity property
\begin{equation}\label{L1weight.contr.Step.1.2}
\int_{\Omega}\big[u(t,x)-v(t,x)\big]\Phi_1(x)\dx\le\int_{\Omega}\big[u(\tau,x)-v(\tau,x)\big]\Phi_1(x)\dx\,,
\qquad\mbox{for all }0\le \tau\le t\,.
\end{equation}
Next, we combine the  weighted smoothing effect \eqref{thm.Upper.PME.Smoothing.2} and the monotonicity inequality \eqref{L1weight.contr.Step.1.2} with the numerical inequality $(u-v)(u^m-v^m)\le m\max\{u^{m-1}\,,\,v^{m-1}\}(u-v)^2$ valid for all $m\ge 1$ and $u,v\ge 0$\,, (recall that $u-v\ge 0$), to get for all $t\ge \tau_0$:
\begin{equation}\label{L1weight.contr.Step.1.3}\begin{split}
&\int_\Omega \big(u^m(t,x)-v^m(t,x)\big)\Phi_1(x)\dx
\le m \max\left\{\|u(t)\|_{\LL^\infty(\Omega)}^{m-1}\,,\,\|v(t)\|_{\LL^\infty(\Omega)}^{m-1}\right\}
    \int_\Omega \big(u(t,x)-v(t,x)\big)\Phi_1(x)\dx\\
&\le \frac{m\, K_4^{m-1}}{(t-\tau_0)^{N(m-1)\vartheta_{1,1}}}
    \max\left\{\|u(\tau_0)\|_{\LL^1_{\Phi_1}(\Omega)}^{2s\vartheta_{1,1}(m-1)}
            \,,\,\|v(\tau_0)\|_{\LL^1_{\Phi_1}(\Omega)}^{2s\vartheta_{1,1}(m-1)}\right\}
    \int_\Omega \big(u(\tau_0,x)-v(\tau_0,x)\big)\Phi_1(x)\dx\\
&:= \frac{\overline{K}[u(\tau_0),v(\tau_0)]}{(t-\tau_0)^{N(m-1)\vartheta_{1,1}}}\int_\Omega \big(u(\tau_0,x)-v(\tau_0,x)\big)\Phi_1(x)\dx\,.
\end{split}\end{equation}
Plugging inequality \eqref{L1weight.contr.Step.1.3} into \eqref{L1weight.contr.Step.1.1} gives for all $t_1,t_0\ge \tau$\,:
\begin{equation}\label{L1weight.contr.Step.1.4}
\begin{split}
&\left|\int_{\Omega}\big[u(t_0,x)-v(t_0,x)\big]\Phi_1(x)\dx - \int_{\Omega}\big[u(t,x)-v(t,x)\big]\Phi_1(x)\dx\right|\\
&\le \lambda_1
\int_{t_0}^{t_1}\frac{\overline{K}[u(\tau_0),v(\tau_0)]}{(t-\tau_0)^{N(m-1)\vartheta_{1,1}}}
\int_\Omega \big(u(\tau_0,x)-v(\tau_0,x)\big)\Phi_1(x)\dx\\
\end{split}
\end{equation}
which implies \eqref{L1weight.contr.estimates.1} using the numerical inequality $|(t-\tau_0)^{2s\vartheta_{1,1}}-(\tau-\tau_0)^{2s\vartheta_{1,1}}|\le |t-\tau|^{2s\vartheta_{1,1}}$\,, valid since $2s\vartheta_{1,1}<1$\,. The constant $K_7=\lambda_1\overline{K}[u(\tau_0),v(\tau_0)]/(2s\vartheta_{1,1})$\,, has the form given in \eqref{L1weight.contr.estimates.1.a}.\qed

\medskip

\noindent{\bf Proof of Theorem \ref{thm.L1weight.exist}.} \\
\noindent$\bullet~$\textsc{Step 1. }In a first step we take bounded and nonnegative initial data $u_0$. In this case we may recover the existence and uniqueness of strong solutions proved in \cite{DPQRV1, DPQRV2} and quoted in Section \ref{sect.statement}.

\noindent$\bullet~$\textsc{Step 2. }We now take general data $u_0\ge 0$\,, and construct a limit solution by approximation from below with $\LL^\infty$-strong solutions. We consider a monotone non-decreasing sequence $0\le u_{0,n}\le u_{0,n+1}\le u_0$\,, with $u_{0,n}\in \LL^\infty(\Omega)$\,, monotonically converging from below to $u_0\in \LL^1_{\Phi_1}(\Omega)$ in the topology of $\LL^1_{\Phi_1}(\Omega)$. By previous results, we know that to every $u_{0,n}$ corresponds a unique strong solution $u_n(t,x)$.
Since $u_n$ is a strong solution, in particular it is a very weak solution in the sense of Definition \ref{Def.Very.Weak.Sol.Dual}\,.  Moreover, since $u_{0,n}\in \LL^\infty(\Omega)$\,, then $u_n(t)\in \LL^\infty(\Omega)$\,, and the comparison holds for these solutions, so that the sequence $u_n(t,x)$ is ordered, namely $u_n(t,x)\le u_{n+1}(t,x)$ for all $x\in \Omega$ and $t>0$\,. We know that $u_n\in \mathcal{S}$\,, therefore all the upper estimates of Theorems \ref{thm.Upper.PME}, \ref{thm.Upper.PME.II} and \ref{thm.Upper.1.PME} hold true for $u_n$\,. \normalcolor Thus, for any fixed $\tau>0$, there exists the limit $u(t,x)=\lim\limits_{n\to \infty} u_n(t,x)$ in $\LL^\infty((\tau,\infty)\times\Omega)$.

Moreover, $u\in C^0([0,\infty)\,:\,\LL^1_{\Phi_1}(\Omega))$. Indeed, using the estimates \eqref{L1weight.contr.estimates.1} for $u_n\le u$, we have that for all $\tau>0$
\begin{equation}\label{Step.1.1.thm.L1weight.exist}\begin{split}
0\le \int_{\Omega}\big[u(\tau,x)-u_n(\tau,x)\big]\Phi_1(x)\dx
&\le \int_{\Omega}\big[u_0(x)-u_{0,n}(t,x)\big]\Phi_1(x)\dx \\
&+ K_7[u_0,u_{0,n}]\,\tau^{2s\vartheta_1}\,\int_{\Omega} \big[u_0(x)-u_{0,n}(x)\big]\Phi_1\dx\xrightarrow[n\to\, \infty]{} 0\,.
\end{split}
\end{equation}

\noindent$\bullet~$\textsc{Step 3. }\textit{The limit solution is a very weak solution. }We have to check that for all $\psi$ such that $\psi/\Phi_1\in C^1_c((0,+\infty): \LL^\infty(\Omega))$ the following identity holds true:
\begin{equation}
\displaystyle \int_0^\infty\int_{\Omega}\AI (u) \,\dfrac{\partial \psi}{\partial t}\,\dx\dt
-\int_0^\infty\int_{\Omega} |u|^{m-1}u\,\psi\,\dx \dt=0.
\end{equation}
Let us fix an admissible test function $\psi$\,. As mentioned above, we know that strong solutions are in particular very weak solutions in the sense of Definition \ref{Def.Very.Weak.Sol.Dual}\,,  so that
\begin{equation}
\displaystyle \int_0^\infty\int_{\Omega}\AI (u_n) \,\dfrac{\partial \psi}{\partial t}\,\dx\dt
=\int_0^\infty\int_{\Omega} u^m_n\,\psi\,\dx \dt\,.
\end{equation}
The proof is concluded once we show that
\begin{equation}\label{Step.2.1.thm.L1weight.exist}
\int_0^\infty\int_{\Omega}\AI (u_n) \,\dfrac{\partial \psi}{\partial t}\,\dx\dt
\xrightarrow[n\to\, \infty]{}
\int_0^\infty\int_{\Omega}\AI (u) \,\dfrac{\partial \psi}{\partial t}\,\dx\dt\,,
\end{equation}
and
\begin{equation}\label{Step.2.2.thm.L1weight.exist}
\int_0^\infty\int_{\Omega} u^m_n\,\psi\,\dx \dt
\xrightarrow[n\to\, \infty]{}
\int_0^\infty\int_{\Omega} u^m\,\psi\,\dx \dt\,.
\end{equation}
\noindent$\bullet~$\textit{Proof of \eqref{Step.2.1.thm.L1weight.exist}. }Recall that $\psi/\Phi_1\in C^1_c((0,+\infty): \LL^\infty(\Omega))$ and say that the time-support is contained in $[t_1,t_2]$\,, and recall that $f\ge 0$ implies $\AI f \ge 0$\,, so that $u\ge u_n$ implies $\AI(u-u_n)\ge 0$ and
\begin{equation*}\begin{split}
& \left|\int_0^\infty\int_{\Omega}\AI (u-u_n)(t,x) \,\partial_t\psi\,\dx\dt\right|
\le \int_{t_1}^{t_2} \left\|\frac{\partial_t\psi}{\Phi_1}\right\|_{\LL^\infty(\Omega)}\int_{\Omega}\Phi_1(x)\AI (u-u_n)(t,x) \,\dx\dt\\
&= \int_{t_1}^{t_2} \left\|\frac{\partial_t\psi}{\Phi_1}\right\|_{\LL^\infty(\Omega)}\int_{\Omega}\AI\Phi_1(x) (u-u_n)(t,x) \,\dx\dt
= \lambda_1\int_{t_1}^{t_2} \left\|\frac{\partial_t\psi}{\Phi_1}\right\|_{\LL^\infty(\Omega)}\int_{\Omega} (u-u_n)(t,x)\Phi_1(x) \,\dx\dt\\
& \le \lambda_1\int_{\Omega} (u_0-u_{0,n})\Phi_1\,\dx\, \int_{t_1}^{t_2} \left\|\frac{\partial_t\psi}{\Phi_1}\right\|_{\LL^\infty(\Omega)}\dt \xrightarrow[n\to\, \infty]{} 0 \\
\end{split}
\end{equation*}
where we have used the fact that $\AI$ is symmetric and in the last step we have used inequality \eqref{Step.1.1.thm.L1weight.exist}\,.

\noindent$\bullet~$\textit{Proof of \eqref{Step.2.2.thm.L1weight.exist}. }Recall that $\psi/\Phi_1\in C^1_c((0,+\infty): \LL^\infty(\Omega))$ and say that the time-support is contained in $[t_1,t_2]\subset (0,\infty)$\,, so that
\begin{equation*}\begin{split}
\left|\int_0^\infty\int_{\Omega} (u^m-u_n^m)(t,x)\,\psi\,\dx \dt\right|
&\le \int_{t_1}^{t_2} \left\|\frac{\psi}{\Phi_1}\right\|_{\LL^\infty(\Omega)}\int_{\Omega}(u^m-u_n^m)(t,x)\,\Phi_1(x) \,\dx\dt\\
&\le m\int_{t_1}^{t_2} \left\|\frac{\psi}{\Phi_1}\right\|_{\LL^\infty(\Omega)} \|u(t)\|_{\LL^\infty(\Omega)}^{m-1}\int_{\Omega}(u-u_n)(t,x)\,\Phi_1(x) \,\dx\dt\\
&\le m\int_{t_1}^{t_2} \left\|\frac{\psi}{\Phi_1}\right\|_{\LL^\infty(\Omega)} \frac{K_1^{m-1}}{t}\int_{\Omega}(u_0-u_{0,n})\,\Phi_1 \,\dx\dt
\xrightarrow[n\to\, \infty]{} 0 \\
\end{split}
\end{equation*}
where we have used the numerical inequality $(u-v)(u^m-v^m)\le m\max\{u^{m-1}\,,\,v^{m-1}\}(u-v)^2$ valid for all $m\ge 1$ and $u,v\ge 0$\,, (recall that $u-u_n\ge 0$)\, and the fact that $\|u_n(t)\|_{\LL^\infty(\Omega)}\le \|u(t)\|_{\LL^\infty(\Omega)}$ for all $t>0$\,, by construction. Finally, we have used the absolute bounds \eqref{thm.Upper.PME.Absolute} as follows: we know that for all $n$ and for all $t>0$ we have $\|u_n(t)\|_{\LL^\infty(\Omega)}\le K_1\,t^{1/(m-1)}$\,, therefore, for all $t>0$, we have that
\begin{equation}\label{final.002}
\|u(t)\|_{\LL^\infty(\Omega)}\le \limsup_n\|u_n(t)\|_{\LL^\infty(\Omega)}\le  K_1\,t^{1/(m-1)}\,.
\end{equation}
Finally, we have also used  inequality \eqref{Step.1.1.thm.L1weight.exist}\,.

\noindent$\bullet~$\textsc{Step 4. }\textit{Uniqueness. }We keep the notations of the previous steps.
Assume that there exist another monotone non-decreasing sequence $0\le v_{0,k}\le v_{0,k+1}\le u_0$\,, with $v_{0,k}\in\LL^\infty(\Omega)$\,, monotonically converging from below to $u_0\in \LL^1_{\Phi_1}(\Omega)$ in $\LL^\infty((\tau,\infty)\times\Omega)$ for all $\tau>0$. By the same considerations as in the proof of Theorem \ref{thm.L1weight.exist} we can show that there exists a solution $v(t,x)\in C^0([0,\infty)\,:\,\LL^1_{\Phi_1}(\Omega))$. We want to show that $u=v$, where $u$ is the solution constructed in the same way from the sequence $u_{0,n}$. We will prove equality by proving that $v\le u$ and then that $u\le v$. To prove that $v\le u$ we use the estimates
\begin{equation}\label{final.001}
\int_{\Omega}\big[v_k(t,x)-u_n(t,x)\big]_+\dx\le \int_{\Omega}\big[v_k(0,x)-u_n(0,x)\big]_+\dx
\end{equation}
which hold for any $u_n(t,\cdot), v_k(t,\cdot)\in \LL^\infty(\Omega)$, by the results of \cite{DPQRV1,DPQRV2}\,. Letting $n\to \infty$ we get that
\[
\lim_{n\to\infty}\int_{\Omega}\big[v_k(t,x)-u_n(t,x)\big]_+\dx
\le \lim_{n\to\infty}\int_{\Omega}\big[v_k(0,x)-u_n(0,x)\big]_+\dx
=\int_{\Omega}\big[v_k(0,x)-u_0(x)\big]_+\dx=0
\]
since $v_k(0,x)\le u_0$ by construction. Therefore also $v_k(t,x)\le u(t,x)$ for $t>0$, so that in the limit $k\to \infty$ we obtain $v(t,x)\le u(t,x)$\,. The inequality $u\le v$ can be obtained simply by switching the roles of $u_n$ and $v_k$\,. The validity of estimates of Proposition $\ref{thm.L1weight.contr}$ is guaranteed by the above limiting process. The comparison holds by taking the limits in inequality \eqref{final.001}\,.

\noindent$\bullet~$\textsc{Step 5. }It remains to show that  the solutions constructed above  belong to the class $\mathcal{S}$, so that the upper and lower bounds of Theorems $\ref{thm.Upper.PME}$, $\ref{thm.Upper.PME.II}$, $\ref{thm.Upper.1.PME}$ and $\ref{thm.Lower.PME}$ apply. We have already shown that $u\in C([0,\infty):\LL^1_{\Phi_1}(\Omega))$ so that we only have to show that  $u(t)\in \LL^p(\Omega)$ for all $t>0$, with $p>N/(2s)$\,. We have already proved that $u(t)\in \LL^\infty(\Omega)$ for all $t>0$, as a consequence of the absolute bounds of Theorem $\ref{thm.Upper.PME}$, see the end of Step 2, formula \eqref{final.002}\,. The scale invariance property follows from uniqueness.  Therefore $u\in \mathcal{S}$.\qed

\section{Appendix}

\subsection{Reminder on Green Functions}\label{SSection.Green}

As already mentioned in the introduction, one of the novelties of this paper is represented by the technique used in the proof of the lower and upper estimates. A main ingredient is the knowledge of good estimates for the Green function. It is known that the Green function of $\A$ satisfies the following estimate for all $x,x_0\in \Omega$:
\begin{equation}\label{Gree.est.0}
G_{\Omega}(x,x_0)\asymp \frac{1}{|x-x_0|^{N-2s}}
\left(\frac{\dist(x,\partial\Omega) }{|x-x_0| }\wedge 1\right)
\left(\frac{\dist(x,\partial\Omega) }{|x-x_0| }\wedge 1\right)\,,
\end{equation}
where $(a\wedge 1)(b\wedge 1)=\min\{1,a,b,ab\}$ for all $a,b\ge 0$\,. Indeed the above estimate can be obtained by the following Heat kernel estimates $H(t,x,y)$ for $s=1$
\begin{equation}\label{Heat.Kernel.est}
H(t,x,x_0)\asymp \left(\frac{\phi_1(x)}{|x-x_0|}\wedge 1\right)
\left(\frac{\phi_1(x_0)}{|x-x_0|}\wedge 1\right)\,\frac{\ee^{-\frac{c_2|x-x_0|^2}{t}}}{t^{N/2}}\,.
\end{equation}
with the help of the formula
\[
G_{\Omega}(x,y)=\int_0^{\infty}\frac{H(t,x,y)}{t^{1-s}}\dt\,.
\]
Recall that in this case $\Phi_1=\phi_1$  is the first eigenfunction of the Dirichlet Laplacian ($s=1$) and satisfies estimates \eqref{Phi1.est}\,, namely $\Phi_1=\phi_1\asymp\dist(\cdot,\partial\Omega)$\,. The upper bounds for the Heat kernel \eqref{Heat.Kernel.est} can be found in \cite{D1, D2, Davies1, Davies2, DS}\,, while the lower bounds have been obtained later in \cite{Zh2002}\,.

It is easy to see that estimates \eqref{Gree.est.0} imply respectively the Type I and Type II estimates \eqref{typeI.Green.est} and \eqref{typeII.Green.est}\,, which have been already stated in Section \ref{ssect.Funct.Setup}.

Notice that here we have focused our attention on the SFL, but  similar estimates for the Green function hold for other classes of operators. This will be investigated  in the forthcoming paper \cite{BV-Paper2}.

We now state the integral estimates which have been used in the proofs of the present paper.

\begin{lem}[Integral Green function estimates I]\label{Lem.Green}Let $G_\Omega$ be the Green function of $\A$. Then, the Type I estimates \eqref{typeI.Green.est} imply that there exist a constant $c_{2,\Omega}(q)>0$ such that
\begin{equation}\label{Lem.Green.est.Upper.I}
\sup_{x_0\in\Omega}\int_{\Omega}G^q_{\Omega}(x , x_0)\dx \le c_{2,\Omega}(q)\qquad\mbox{for all $0<q<\dfrac{N}{N-2s}$\,,}
\end{equation}
Moreover, the Type II estimates \eqref{typeII.Green.est} imply that there exist constants $c_{3,\Omega}(q),c_{4,\Omega}(q)>0$ such that for all $0< q <\frac{N}{N-2s}$
\begin{equation}\label{Lem.Green.est.Upper.II}
c_{3,\Omega}(q) \Phi_1(x_0)\le \left(\int_{\Omega}G^q_{\Omega}(x , x_0)\dx\right)^{\frac{1}{q}} \le   c_{4,\Omega}(q) \B_q(\Phi_1(x_0))\,,
\end{equation}
where the function $\B_q:[0,\infty)\to[0,\infty)$ is defined as follows:
\begin{equation}\label{Lem.Green.est.Upper.II.B}
\B_q(\Phi_1(x_0))\left\{\begin{array}{lll}
\Phi_1(x_0)\,, & \qquad\mbox{for any }0< q <\frac{N}{N-2s+1}\,,\\[2mm]
\Phi_1(x_0)\, \big|\log\Phi_1(x_0)\big|^{\frac{1}{q}}\,, & \qquad\mbox{for }q = \frac{N}{N-2s+1}\,,\\[2mm]
\Phi_1(x_0)^{\frac{d-q(d-2s)}{q }}\,, & \qquad\mbox{for any }\frac{N}{N-2s+1}<q<\frac{N}{N-2s}\,.\\
\end{array}\right.
\end{equation}
Moreover, for all $f\in \LL^1_{\Phi_1}(\Omega)$\,, the Type II estimates \eqref{typeII.Green.est} imply for all $x_0\in \Omega$
\begin{equation}\label{Lem.Green.est.Lower.II}
\int_\Omega f(x)G_{\Omega}(x,x_0)\dx\ge c_{0,\Omega} \Phi_1(x_0) \|f\|_{\LL^1_{\Phi_1}(\Omega)}\,,
\end{equation}
The constants $c_{i,\Omega}(\cdot)$\,, $i=2,3,4$\,, may depend on $s,N, q$ and have an explicit expression given at the end of the proof.
\end{lem}

The proof of this result is  technical and long, but not difficult. We will prove a more general version of the above Lemma in \cite{BV-Paper2}.

As a consequence of the above estimates, for any $f\in \LL^\infty(\Omega)$ we have that
\begin{equation}\label{rem1.1.Green}
c_{0,\Omega}\Phi_1(x_0) \|f\|_{\LL^1_{\Phi_1}(\Omega)}\le \int_\Omega f(x)G_{\Omega}(x,x_0)\dx\le c_{4,\Omega}\|f\|_{\LL^\infty(\Omega)}\B_1(\Phi_1(x_0))
\end{equation}
with
\begin{equation}\label{rem1.2.Green}
\B_1(\Phi_1(x_0))=\left\{\begin{array}{lll}
\Phi_1(x_0)\,, & \qquad\mbox{for any } 2s>1\,,\\[2mm]
\Phi_1(x_0)\, \big|\log\Phi_1(x_0)\big| \,, & \qquad\mbox{for }2s=1\,,\\[2mm]
\Phi_1(x_0)^{2s}\,, & \qquad\mbox{for any }2s<1\,.\\
\end{array}\right.
\end{equation}
It is clear at this point that the upper bounds are sharp when $2s>1$\,, since the powers of $\Phi_1$ match. When $2s\le 1$\,, the situation changes, and the above boundary behaviour is not sharp. We have to work a bit more to obtain matching powers of $\Phi_1$. We remark that as a consequence of the above estimates, when $2s\le 1$ we have
\begin{equation}\label{rem1.3.Green}
\int_\Omega f(x)G_{\Omega}(x,x_0)\dx\le  c_{4,\Omega}\|f\|_{\LL^\infty(\Omega)}\left\{\begin{array}{lll}
\frac{1}{\varepsilon}\Phi_1^{1-\varepsilon}(x_0)\,, & \qquad\mbox{for $2s=1$ and for all $\varepsilon>0$}\,,\\[2mm]
\Phi_1(x_0)^{2s}\,, & \qquad\mbox{for any }2s<1\,.\\
\end{array}\right.
\end{equation}
\begin{lem}[Integral Green function estimates II]\label{Lem.Green.2}Let $m>1$ and $G_\Omega$ be the Green function of $\A$. If for all $x_0\in\Omega$
\begin{equation}\label{Lem.Green.2.hyp}
u^m(x_0)\le \kappa_0\int_{\Omega} u(x)G_\Omega(x,x_0)\dx\,.
\end{equation}
then, the Type II estimates \eqref{typeII.Green.est} imply that there exist a constant $c_{5,\Omega}>0$ such that for all $x_0\in \Omega$
\begin{equation}\label{Lem.Green.est.Upper.II.b}
u^m(x_0)\le \kappa_0\int_{\Omega} u(x)G_\Omega(x,x_0)\dx\le c_{5,\Omega}^m\kappa_0^{\frac{m}{m-1}}\Phi_1(x_0)\,.
\end{equation}
The constant $c_{5,\Omega} $\,, may depend on $s,N,m$ but not on $u$ nor on $\kappa_0$, and have an explicit expression given in the  proof.
\end{lem}
\noindent {\sl Proof.~}We split several steps.

\noindent$\bullet~$\textsc{Step 1. }\textit{First boundary estimates. }An immediate consequence of hypothesis \eqref{Lem.Green.2.hyp} is the following absolute upper bound:
\[
u^m(x_0)\le \kappa_0\int_{\Omega} u(x)G_\Omega(x,x_0)\dx \le \kappa_0\,c_{4,\Omega}\|u\|_{\LL^\infty(\Omega)}
\]
where we have used bounds \eqref{Lem.Green.est.Upper.I} for $q=1$\,, namely that $\sup\limits_{x_0\in\Omega}\|G_\Omega(x_0,\cdot)\|_{\LL^1(\Omega)}\le c_{4,\Omega}$\,. Taking the supremum in $x_0$ in the above expression, gives the absolute bound $\|u\|_{\LL^\infty(\Omega)}\le (\kappa_0\,c_{4,\Omega})^{\frac{1}{m-1}}$\,. Finally, using the latter absolute bound, together with hypothesis \eqref{Lem.Green.2.hyp} and inequality \eqref{Lem.Green.est.Upper.II} with $q=1$, we obtain
\begin{equation}\label{Step.0.2.Lem.Green.2}\begin{split}
u^m(x_0)
&\le \kappa_0\int_{\Omega} u(x)G_\Omega(x,x_0)\dx
 \le \kappa_0 c_{4,\Omega}\|u\|_{\LL^\infty(\Omega)} \B_1(\Phi_1(x_0)) 
\end{split}
\end{equation}
Therefore, using \eqref{rem1.3.Green}, we have proven that for all $x_0\in \Omega$\,, we have
\begin{equation}\label{Step.0.3.Lem.Green.2}
u(x_0)\le \kappa_1\Phi_1^{\nu_1}(x_0)\,,
\end{equation}
with
\begin{equation}\label{Step.0.4.Lem.Green.2}
\kappa_1:=\tilde{\kappa}_1\,\kappa_0^{\frac{1}{m-1}}:=(\kappa_0\,c_{4,\Omega})^{\frac{1}{m-1}}\quad\mbox{and}\quad
\nu_1:=\left\{\begin{array}{lll}
1\,, & \qquad\mbox{for any } 2s>1\,,\\[2mm]
\frac{1-\varepsilon}{m}\,, & \quad\mbox{for $2s=1$ and for all $\varepsilon\in (0,1]$}\,,\\[2mm]
\frac{2s}{m}\,, & \quad\mbox{for any }2s<1\,.\\
\end{array}\right.
\end{equation}
therefore inequality \eqref{Lem.Green.est.Upper.II} is already proved in the case $2s>1$\,. From now on we concentrate in the case $2s\le 1$\,. We will split two cases, namely $2s<1$ and $2s=1$\,.

\noindent$\bullet~$\textsc{Step 2. }We first deal with the case $2s<1$. Combining inequality \eqref{Step.0.3.Lem.Green.2} with hypothesis \eqref{Lem.Green.2.hyp} gives
\[
\begin{split}
&u^m(x_0)\le \kappa_0 \int_{\Omega}  u(x) G_\Omega(x,x_0)\dx
\le\kappa_0\kappa_1\int_{\Omega}  \Phi_1^{\nu_1} G_\Omega(x,x_0)\dx \\
&\le\kappa_0^{\frac{m}{m-1}}c_{4,\Omega}^{\frac{1}{m-1}}
    \left(\int_{\Omega}  G_\Omega(x,x_0)\dx \right)^{1-\nu_1}\left(\int_{\Omega}  \Phi_1 G_\Omega(x,x_0)\dx \right)^{\nu_1}
\le \kappa_0^{\frac{m}{m-1}}c_{4,\Omega}^{\frac{1}{m-1}+1-\nu_1}  \Phi_1(x_0)^{2s(1-\nu_1)+\nu_1}
\end{split}
\]
since we recall that since $2s<1$ estimates \eqref{rem1.1.Green} give $\int_{\Omega}  G_\Omega(x,x_0)\dx \le c_{4,\Omega}\Phi_1(x_0)^{2s}$ and we also recall that  $\Phi_1(x_0)=\int_{\Omega}  \Phi_1 G_\Omega(x,x_0)\dx$\,.
\begin{equation}\label{Step.2b.4.Lem.Green.2}
u(x_0)\le \kappa_0^{\frac{1}{m-1}}\widetilde{\kappa}_2\Phi_1(x_0)^{\nu_2}
\end{equation}
with $\nu_2=\frac{2s(1-\nu_1)+\nu_1}{m}=\nu_1(1+\frac{1-2s}{m})>\nu_1$ and $\widetilde{\kappa}_2=c_{4,\Omega}^{\frac{1}{m(m-1)}+\frac{1-\nu_1}{m}}$\,. Iterating the above process $n$ times, gives
\begin{equation}\label{Step.2b.4.Lem.Green.2.b}
u(x_0)\le \kappa_0^{\frac{1}{m-1}}\widetilde{\kappa}_n\Phi_1(x_0)^{\nu_n}\qquad\mbox{with}\qquad
\nu_n=\nu_1\left(1+\frac{1-2s}{m}\right)^n\wedge 1
\end{equation}
where $\kappa_n\le c_{4,\Omega}^{\sigma_n}$ for some finite $\sigma_n$\,.  We have taken $n$ to be the biggest integer such that $\nu_{n+1}>1$\,, namely
\[
n+1>\frac{\log{\frac{1}{\nu_1}}}{\log\left(1+\frac{1-2s}{m}\right)}=\frac{\log{\frac{m}{2s}}}{\log\left(1+\frac{1-2s}{m}\right)}
\]
so that $\nu_n=1$\,, therefore \eqref{Step.2b.4.Lem.Green.2.b} gives inequality \eqref{Lem.Green.est.Upper.II} with $c_{5,\Omega}=c_{4,\Omega}^{\sigma_n}$ for some finite $\sigma_n$\,.

\noindent$\bullet~$\textsc{Step 3. }We deal with the case $2s=1$. Let us fix $x_0\in \Omega$, and $R_0=\Phi_1(x_0)\le \overline{R}=\|\Phi_1\|_{\LL^\infty(\Omega)}+\diam(\Omega)\,,$ so that for any $x_0\in\Omega$ we have $\Omega\subseteq B_{\overline{R}}(x_0)$\,. Notice that it is not restrictive to assume $0\le \Phi_1(x_0)\le 1$ since we know by estimates \eqref{Phi1.est} that $\Phi_1(x)\asymp \big(\dist(x, \partial\Omega)\wedge 1\big)$ for all $x\in\Omega$ and we are interested in the boundary behaviour, i.e. to the values $\dist(x, \partial\Omega)<<1$; when $\Phi_1(x_0)\ge 1$, i.e. far from the boundary, already have estimates \eqref{Lem.Green.est.Upper.I}. Recall now the upper part of Type II estimates \eqref{typeII.Green.est}, that can be rewritten in the equivalent form
\begin{equation}\label{Step.1.1.Lem.Green.2}
G_{\Omega}(x,x_0)\le
\frac{c_{1,\Omega}}{|x-x_0|^{d-2s}}
\left\{\begin{array}{cl}
\dfrac{\Phi_1(x_0)}{|x-x_0|} &\quad\mbox{for any }x\in \Omega\setminus B_{R_0}(x_0)=B^c_{R_0}(x_0)\\
1&\quad\mbox{for any }x\in B_{R_0}(x_0)\\
\end{array}\right.
\end{equation}
Next, we recall that $\Phi_1\in C^1(\overline{\Omega})$\,, so that there exists $k_1>0$ so that for all $x\in\Omega$
\begin{equation}\label{Step.1.2.Lem.Green.2}
\Phi_1(x)\le \Phi_1(x_0)+k_1|x-x_0| \le (1+k_\gamma) \left\{\begin{array}{cl}
|x-x_0|   &\quad\mbox{for any }x\in \Omega\setminus B_{R_0}(x_0)=B^c_{R_0}(x_0)\\
\Phi_1(x_0)&\quad\mbox{for any }x\in B_{R_0}(x_0)\\
\end{array}\right.\,.
\end{equation}
Joining the above estimates we obtain
\begin{equation}\label{Step.1.3.Lem.Green.2}
\Phi_1(x)^{\nu_1} G_{\Omega}(x,x_0)\le
\frac{c_{1,\Omega}(1+k_\gamma)}{|x-x_0|^{d-2s}}
\left\{\begin{array}{cl}
\dfrac{\Phi_1(x_0)}{|x-x_0|^{1-\nu_1}} &\quad\mbox{for any }x\in \Omega\setminus B_{R_0}(x_0)=B^c_{R_0}(x_0)\\
\Phi_1(x_0)&\quad\mbox{for any }x\in B_{R_0}(x_0)\\
\end{array}\right.
\end{equation}
We now recall that $R_0=\Phi_1(x_0)$ and that $2s=1$: we use \eqref{Step.1.3.Lem.Green.2} to estimate
\[\begin{split}
\int_{\Omega}  u(x) & G_\Omega(x,x_0)\dx
\le \kappa_1 \int_{B_{R_0}(x_0)} \Phi_1^{\nu_1}(x)G_\Omega(x,x_0)\dx+\kappa_1 \int_{B^c_{R_0}(x_0)} \Phi_1^{\nu_1}(x)G_\Omega(x,x_0)\dx\\
&\le \kappa_1 c_{1,\Omega}(1+k_1)\left[\int_{B_{R_0}(x_0)} \frac{\Phi_1^{\nu_1}(x_0)\dx}{|x-x_0|^{d-2s}}
+ \int_{B_{\overline{R}}(x_0)\setminus B_{R_0}(x_0)} \dfrac{\Phi_1(x_0)\dx}{|x-x_0|^{d-2s+ (1-\nu_1)}}\right]\\
&= \kappa_1 c_{1,\Omega}(1+k_1)\omega_N\left[\frac{R_0^{2s}}{2s}\Phi_1^{\nu_1}(x_0)
            + \Phi_1(x_0)\frac{\overline{R}^{\nu_1}-R_0^{\nu_1}}{\nu_1}\right]\le \kappa_1 c_{1,\Omega}(1+k_1)\omega_N\left[1+\frac{1}{\nu_1}\right]\,\overline{R}^{\nu_1}\,\Phi_1(x_0)\\
\end{split}
\]
Combining the above inequality with hypothesis \eqref{Lem.Green.2.hyp} gives \eqref{Lem.Green.est.Upper.II} with
$c_{5,\Omega}=c_{4,\Omega}^{\frac{1}{m-1}} c_{1,\Omega}(1+k_1)\omega_N\left(1+\frac{1}{\nu_1}\right)\,(\|\Phi_1\|_{\LL^\infty(\Omega)}+\diam(\Omega))^{\nu_1}$\,.\qed

\noindent\textbf{Remark. }An important consequence of this Lemma is that it makes possible to understand the sharp boundary behaviour of solutions to the Dirichlet elliptic problem $\A u^m= \lambda u$\,, cf. Section \ref{sect.Elliptic.Prob}.\normalcolor
\normalcolor


\section{Comments, extensions and open problems}\label{sect.Gener}

\noindent $\bullet$ We want to stress  an interesting conclusion of our results:  the boundary behaviour of our solutions is dictated by the first eigenfunction $\Phi_1$ of operator $\cal L$. When dealing with other linear operators, or with less regular domains, $\Phi_1$ may have a different boundary behaviour, namely $\big(\dist(\cdot, \partial\Omega)^\gamma\wedge 1\big)$, but the main results can be shown to hold, with appropriate expressions.  Roughly speaking, $u(t,\cdot)$ behaves like the distance to the boundary at a the power $\gamma/m$. The detailed analysis of this extension will appear in a forthcoming paper \cite{BV-Paper2}.\normalcolor

\medskip

\noindent $\bullet$  When the equation is posed in the whole space $\Omega=\RR^N$, the questions  we have discussed here have been treated in the recent literature, cf. \cite{BV2012, DPQRV1, DPQRV2}. While a number of general ideas are similar, the details and techniques are quite different. The difficulties found in applying previous ideas have led to the new methods used here.

\medskip

\noindent $\bullet$ {\sc The case  $m=1$}. In this limit we obtain the fractional linear heat equation. The theory is simpler since the solution can be obtained by means of the representation formula in terms of the heat kernel $H_s(t,x,y)$ or in terms of the eigenvalues and eigenfunctions $(\lambda_k^s,\Phi_k)$: recall that $\widehat{u}_{k}=\int_\Omega u_0\Phi_k\dx$, so that
\[
u(t,x)=\int_\Omega H_s(t,x,y) u_0(y)\dy = \sum_{k=1}^\infty \ee^{-\lambda_k^s t}\, \widehat{u}_{k}\,\Phi_k(x) \,,
\]
from which it follows by standard methods that $u(t,x)\asymp \ee^{-\lambda_1^s t} \Phi_1(x)$ for all large times\,, cf. \cite{Davies1, Davies2}\,. Notice that there cannot be absolute upper bounds since the equation is linear. There exists a linear version of the positivity result of Theorem \ref{thm.Lower.PME}\,.

\medskip

\noindent $\bullet$ {\sc The case $s=1$}. In this limit case we obtain the standard version of the  porous medium equation.  Our proofs above are stable under this limit. Some of these results are new, like elliptic Harnack inequalities, or backward smoothing effects.

\medskip

\noindent $\bullet$ {\sc Solutions with any sign}. Due to the property of comparison, extending the upper bounds  to solutions of any sign is easy. Here is the argument: recall first that if $u$ is a solution\,, also $-u$ is a solution. Then, consider the nonnegative solution $u^+$ corresponding to $u_0^+=\max\{u_0,0\}$. Thus by comparison, $u\le u^+$\,, since $u_0\le u_0^+$. Consider also the nonnegative solution $u_-$ corresponding to $u_0^-=-\min\{u_0,0\}$. Then by comparison, $-u\le u^-$\,, since $-u_0\le -u_0^-$, and we get  $-u^-\le u\le u^+$, so that $|u|\le \max\{u^+\,,\,u^-\}$\,.

\medskip

\noindent $\bullet$ {\sc More general operators.} Since the derivation of our estimates proceeds on a rather general level, (a large part of) the conclusions should apply to more general operators, like variations of our spectral laplacian operator, such as the case fractional powers of operators with measurable coefficients (under appropriate conditions), or fractional powers of the Laplace-Beltrami Laplacian on Riemannian manifolds.  We will address this question in the forthcoming paper \cite{BV-Paper2}.\\[3mm]

{\sc {\large Acknowledgment.}} Both authors partially funded by Project MTM2011-24696 (Spain).

\end{document}